\documentclass{article}

\usepackage{amssymb}
\usepackage{graphicx}
\usepackage{amsmath}
\usepackage{theorem}

\newtheorem{theorem}{Theorem}

\newtheorem{corollary}[theorem]{Corollary}

\newtheorem{definition}[theorem]{Definition}

\newtheorem{lemma}[theorem]{Lemma}

\newtheorem{proposition}[theorem]{Proposition}
\newtheorem{remark}[theorem]{Remark}

\newenvironment{proof}[1][Proof]{\textbf{#1.} }{\ \rule{0.5em}{0.5em}}

\begin{document}

\title{Typical support and Sanov large deviations of correlated states}
\author{Igor Bjelakovi\'{c}$^{2,3}$ \and Jean-Dominique Deuschel$^{2}$ \and %
Tyll Kr\"{u}ger$^{1,2,4}$ \and Ruedi Seiler$^{2}$ \and Rainer
Siegmund-Schultze$^{1,2}$ \and Arleta Szko\l a$^{1,2}$ \\
$^{1}${\footnotesize Max Planck Institute for Mathematics in the Sciences}\\
{\footnotesize Inselstrasse 22, 04103 Leipzig, Germany}\\
$^{2}${\footnotesize Technische Universit\"{a}t Berlin} \\
{\footnotesize Fakult\"at II - Mathematik und Naturwissenschaften} \\
{\footnotesize Institut f\"ur Mathematik MA 7-2} \\
{\footnotesize Stra\ss e des 17. Juni 136,} {\footnotesize 10623 Berlin,
Germany}\\
$^{3}${\footnotesize Heinrich-Hertz-Chair for Mobile Communication}\\
{\footnotesize Technische Universit\"{a}t Berlin }\\
{\footnotesize Werner-von-Siemens-Bau (HFT 6) }\\
{\footnotesize Einsteinufer 25, 10587 Berlin, Germany}\\
$^{4}${\footnotesize Universit\"{a}t Bielefeld}\\
{\footnotesize Fakult\"{a}t f\"{u}r Physik}\\
{\footnotesize Universit\"{a}tsstr. 25, 33619 Bielefeld, Germany}}
\maketitle

\begin{abstract}
Discrete stationary classical processes as well as quantum lattice states
are asymptotically confined to their respective typical support, the
exponential growth rate of which is given by the (maximal ergodic) entropy.
In the iid case the distinguishability of typical supports can be
asymptotically specified by means of the relative entropy, according to
Sanov's theorem. We give an extension to the correlated case, referring to
the newly introduced class of HP-states.
\end{abstract}


\section{Introduction}

A relevant notion on the interface of classical discrete probability theory
and information theory is that of typical subsets. For the quantum
extensions of these fields there is a corresponding notion: typical
subspaces.

The general picture is that a stationary process (\emph{state} in the case
of quantum lattice systems) is asymptotically -i.e. observing a large finite
interval- more and more confined to its typical support. The size of this
support has an exponential growth rate (possibly zero) given by the
essential supremum of the \emph{entropies} of the ergodic components. In the
classical situation this is the content of the Shannon-McMillan theorem. It
clarifies the importance of \emph{Shannon entropy} for several fields, from
data transmission and compression to statistical mechanics or complexity
theory.

Under the much stronger condition of complete independence \emph{Sanov's
theorem} (see \cite{orgsanov} or \cite{covth})\ specifies the exponential
rate of this confinement of a classical iid process to its own typical set,
or equivalently, the rate of avoidance of the supports of \emph{all other}
iid processes' typical sets. This large deviations result is usually seen as
a result on \emph{empirical distributions}, as in its formulation a
particular instance of typical set appears: typical for an iid process are
realizations with an empirical distribution close to the probability
distribution underlying this very process, see ch. 3.2 in Deuschel and
Stroock \cite{deustr}.

In the iid case Sanov's theorem significantly extends the assertion of the
Shannon-McMillan theorem. In fact, taking the equidistribution as \emph{%
reference measure}, it follows from Sanov's theorem that there is a \emph{%
universal typical set} sequence of approximate size $e^{nh}$ for \emph{all }%
iid processes with (base $e$) entropy less than $h$. It is well-known in the
classical situation that this extends to the general ergodic case (since
there exist universal compression schemes like the Lempel-Ziv algorithm).
This universality result was generalized to the quantum case by Kaltchenko
and Yang \cite{ky}, using a nice 'rotation technique' and the quantum
Shannon-McMillan theorem \cite{q-sm}.

From the point of view of \emph{statistical hypothesis testing} Sanov's
theorem asserts that there is a universally typical set sequence for any set
of iid probability distributions (null hypothesis), separating it optimally
from \emph{any} other set of iid processes (alternative hypothesis) at a
rate arbitrarily close to the infimum of the \emph{relative entropies}
between probability measures from the two hypotheses. So in the classical
case Sanov's theorem expresses a twofold universality in the choice of the
typical sets.

The special case of Sanov's theorem with both hypotheses consisting of only
one probability distribution each, is usually called \emph{Stein's lemma}$.$

As already emphasized in \cite{sanov-i.i.d}, when passing from the classical
to the quantum case, the universality mentioned above gets partially lost:
there exists no longer a sequence of typical subspaces (of the underlying
finite dimensional Hilbert spaces for the $n$-blocks of the system), which
would work universally, whatever the reference states are. Consequently,
speaking in the hypothesis testing terminology, for the alternative
hypothesis only one process/state is admitted here. Universality with
respect to the null hypothesis states is maintained, however. Also, in the
quantum situation it is no longer possible to originate Sanov's theorem on
the concept of empirical distributions (states), see \cite{sanov-i.i.d},
chapter 4.

We mention here that the main techniques needed to generalize Sanov's
theorem to the \emph{iid} quantum case were already presented in Hayashi 
\cite{hayashi} (1997), and in Hayashi \cite{hayashi2} (2002) an equivalent
result is shown. The authors of the present paper regrettably were not aware
of this part of Hayashi's work during the preparation of \cite{sanov-i.i.d}.

It is the aim of this paper as a continuation of \cite{sanov-i.i.d} to
extend the assertion of Sanov's theorem in several directions. This concerns
the classical case, too, but the main focus is on the quantum situation.

First, the restriction to the uncorrelated case is substantially alleviated: 
\emph{No condition besides stationarity} is imposed to the processes/states $%
P$\ of the null hypothesis. As for the alternative hypothesis (reference
measure/state) $Q$, even stationarity is not assumed. The only requirements
are the existence of relative entropy rates $h(W,Q)\leq +\infty $ for the
ergodic components $W$ occurring in the null hypothesis set, and the
validity of the upper bound (\emph{achievability part}) in Stein's lemma
concerning $W$ and $Q$ (see Theorems \ref{stat-sanov} resp. \ref%
{stat-sanov-class} in the classical situation). These are, in a sense, \emph{%
minimal} requirements, since Stein's lemma is a trivial consequence of
Sanov's theorem obtained by forgetting about universality.

As an application of this general result we consider the case that a certain
(admittedly very strong) mixing condition holds for the reference process $Q$%
. Observe that the very existence of the relative entropy rate for
correlated processes can only be guaranteed \emph{in terms of mixing
conditions}, if the reference process is particulary strong mixing. Shields 
\cite{shields} gives an example where the reference process is even \emph{%
maximally mixing} in the sense of Dynamical Systems theory ($B$-process,
i.e. isomorphic to an iid process), but nonetheless there exists no
asymptotic rate of the relative entropy. Though the mixing condition upon
the reference processes is very strong ($\ast $-mixing, cf. \cite{bhk}, or 
\cite{bradley}, where it is called $\psi $-mixing), the class of aperiodic
irreducible \emph{Markov processes} on a finite state space is covered. In
this Markov case aperiodicity is necessary and sufficient for mixing, but
not needed for Sanov's theorem, showing that $\ast $-mixing is far from
being a necessary condition for a Sanov type theorem. In fact, in the
classical case a kind of average-mixing would yield the result, cf.
condition ($\widehat{\text{U}}$) on page 86 of \cite{deustr}.

We also mention, that in the classical case a usual mixing condition to
derive large deviation results is \emph{hypermixing}, cf. chapter 5.4 in 
\cite{deustr}.

Secondly, we generalize the classical Sanov's theorem to the (correlated) 
\emph{quantum situation }(in Hayashi \cite{hayashi2} and later in \cite%
{sanov-i.i.d} the quantum iid case was considered). In fact, since the
classical assertion is a special case of the quantum theorem, we only prove
the latter. Again, the reference state only needs to fulfil the two minimal
conditions mentioned above. We refer to those as HP-condition. The states
forming the null hypothesis have to be stationary only.

It would be interesting to specify the set of all states which fulfil the
HP-condition with respect to \emph{any} ergodic (null hypothesis) state. We
call these states \emph{HP-states}. As already said, this set comprises all $%
\ast $-mixing states, but can be expected to be much larger.

In the classical situation we remind the reader of an interesting example by
Xu \cite{xuexample}: There exists a $B$-process (i.e. again maximally mixing
in the sense of Dynamical Systems) $Q$ which has the property that the
relative entropy rate $h(P,Q)$ exists and \emph{is zero }for \emph{any}
stationary process $P$. So this process cannot be separated at exponential
speed from an arbitrary other stationary process.

It would be interesting to find conditions weaker than $\ast $-mixing
ensuring exponential separability in the case that the relative entropy rate
is positive.

In the presented form, the quantum Sanov's theorem comprises and extends
several earlier results on typical subspaces and their connection with the
von Neumann entropy and relative entropy.

In particular, the result \cite{sanov-i.i.d} of the present authors (which
was preceeded by Hayashi \cite{hayashi2}) is extended to the correlated
case. The quantum Shannon-McMillan theorem \cite{q-sm} is covered and
extended from the ergodic to the general stationary situation. The
universality result of Kaltchenko and Yang \cite{ky} is covered, too, by
using the tracial state as reference state in Theorem \ref{stat-sanov}. In
fact this Kaltchenko-Yang universality is a main ingredient in our proof.
The quantum Stein's lemma (see \cite{ogawanagaoka}, \cite{petzbook}, chapter
1.1 for the iid case, \cite{igrai} for the case of ergodic null hypothesis
states) is covered and extended to the case of correlated reference states.
Results of Hiai and Petz \cite{hiai-petz1}, \cite{hiai-petz2} are completed
in the sense that their bound is shown to be sharp, which means that it is
asymptotically optimal, and the condition of \emph{complete ergodicity}
concerning the null hypothesis is dropped. In particular, the case of
irreducible aperiodic \emph{algebraic }(reference)\emph{\ states }on a
quasi-local algebra over a finite-dimensional $C^{\ast }$-algebra $\mathcal{A%
}$ (also called finitely correlated states) considered in \cite{hiai-petz2}
is covered by $\ast $-mixing. We mention, that Hiai and Petz emphasize in 
\cite{hiai-petz2} that they derive almost all assertions using $\ast $%
-mixing, only.

As already emphasized, the quantum Sanov' theorem is a special type of
quantum large deviations result. We refer the reader to some other work in
this direction, see Lebowitz, Lenci and Spohn \cite{leblensp}, Lenci and
Rey-Bellet \cite{lenluc}, Neto\v{c}n\'{y} and Redig \cite{netred}, De Roeck,
Maes and Neto\v{c}n\'{y} \cite{rmn}.

\bigskip

We give a short account of the principal steps to prove the main result.

In chapter \ref{chsteinsan} we show that 'one half' of Stein's lemma, namely
the assumed achievability of the relative entropy rate as separation rate,
already implies Sanov's theorem.

First it is shown that the optimality of the relative entropy rate (as
separation rate) is a consequence of its achievability. In fact, for two
states $\Psi ,\Phi $ such that $s(\Psi ,\Phi )$ and $s(\Psi )$ exist, the
quantity $-n(s(\Psi ,\Phi )+s(\Psi ))$ is the asymptotic average of the
logarithmic eigenvalues of $D_{\Phi ^{(n)}}$, which denotes the density
operator of the local state $\Phi ^{(n)}$ on the disrete interval of length $%
n$, with respect to the probability measure generated on the corresponding
eigenvectors by the operator $D_{\Psi ^{(n)}}$. On the other hand, the
achievability part of Stein's lemma implies that $-n(s(\Psi ,\Phi )+s(\Psi
)) $ is also an essential \emph{upper} bound for these logarithmic
eigenvalues. The key tool to show the latter is Lemma \ref{mainest}. Now,
roughly speaking, with the asymptotic average being the asymptotic upper
bound, it must be an asymptotic lower bound, too. This observation yields a 
\emph{relative AEP }(asymptotic equipartition property) for the logarithmic
eigenvalues of $D_{\Phi ^{(n)}}$: the vast majority of them (with respect to
the considered probability distribution) is close to $-n(s(\Psi ,\Phi
)+s(\Psi ))$. Because for ergodic $\Psi $ by the quantum Shannon-McMillan
theorem the relevant dimension of the corresponding subspace of eigenvectors
of $D_{\Phi ^{(n)}}$ is close to $e^{ns(\Psi )}$, it easily follows now
(applying Lemma \ref{mainest} once again) that the optimally separating
subspaces can essentially be described as those which are close to the span
of the mentioned eigenvectors of $D_{\Phi ^{(n)}}$ fulfilling the relative
AEP: the $\Phi ^{(n)}$-expectation is close to $e^{ns(\Psi )}\cdot
e^{-n(s(\Psi ,\Phi )+s(\Psi ))}=e^{-ns(\Psi ,\Phi )}$.

Next we make use of the proven relative AEP, combined with Kaltchenko and
Yangs universality result to show Sanov's theorem: We subdivide the null
hypothesis set into small slices of almost constant value of the 'mixed'
term (aka cross-entropy) $s_{\text{mix}}:=s(\Psi ,\Phi )+s(\Psi )=-\lim 
\frac{1}{n}$Tr$D_{\Psi ^{(n)}}\log D_{\Phi ^{(n)}}$, and within these slices
the entropy rate is bounded from above by $s_{\text{mix}}-\inf_{\Psi }s(\Psi
,\Phi )$. Then, by Kaltchenko-Yang universality, there exists a common
support of dimension $\approx e^{n(s_{\text{mix}}-\inf_{\Psi }s(\Psi ,\Phi
))}$ which by the relative AEP can be chosen to consist of eigenvectors of $%
D_{\Phi ^{(n)}}$ with eigenvalues close to $e^{-ns_{\text{mix}}}$. So this
common support has an asymptotic $\Phi ^{(n)}$-expectation close to $e^{n(s_{%
\text{mix}}-\inf_{\Psi }s(\Psi ,\Phi ))}\cdot e^{-ns_{\text{mix}%
}}=e^{-n\inf_{\Psi }s(\Psi ,\Phi )}$.

This essentially proves Sanov's theorem under the HP-condition.

In chapter \ref{starmix} we prove that $\ast $-mixing implies the
HP-condition, hence the quantum Sanov's theorem. The idea is borrowed from 
\cite{hiai-petz2}:\ under $\ast $-mixing, the reference state $\Phi $ is
sufficiently close to some block-iid state, so that we may apply the
techniques developed in \cite{hiai-petz1} and \cite{igrai} in order to prove
the achievability part of Stein's lemma.

In chapter \ref{chstat} we use the ergodic decomposition of stationary
states to extend our results to the case where the null hypothesis states
are only assumed stationary.

\section{Basic settings and notations}

As announced in the introduction, we address both the classical and the
quantum situation. Let us first consider the classical case. Let a finite
set $A$ of symbols be given. We deal with processes $P$ on $[A^{\mathbb{Z}},%
\mathfrak{A}^{\mathbb{Z}}]$, where $\mathfrak{A}^{\mathbb{Z}}$ denotes the $%
\sigma $-field over $A^{\mathbb{Z}}$ which is generated by finite
dimensional cylinders. We denote the set of all processes by $\mathcal{P}(A^{%
\mathbb{Z}})$. Let $P^{(n)}$ denote the marginal of a process $P$,
restricted to the positive (time) indices $\{0,1,...,n-1\}\subset \mathbb{Z}$%
.

The relative entropy rate between two processes $P,Q$ is defined as 
\begin{equation*}
h(P,Q):=\lim_{n\rightarrow \infty }\frac{1}{n}H(P^{(n)},Q^{(n)})
\end{equation*}
whenever this limit exists in $\overline{\mathbb{R}}_{+}=\mathbb{R}_{+}\cup
\{+\infty \}$. Here $H(\cdot ,\cdot )$ denotes the relative entropy of two
probability measures given on a finite set.

If $Q\in \mathcal{P}(A^{\mathbb{Z}}),\ \Omega \subseteq \mathcal{P}(A^{%
\mathbb{Z}})$ and $h(P,Q)$ exists for each $P\in \Omega $ we write $h(\Omega
,Q)$ for $\inf_{P\in \Omega }h(P,Q)$.

The following very strong mixing property of $Q$ was introduced by Blum,
Hanson and Koopmans \cite{bhk} (referred to as $\psi $-mixing in the survey
paper \cite{bradley}), which implies the existence of the relative entropy
rate $h(P,Q)$ for any stationary $P$ (see \cite{hiai-petz2}, where the more
general quantum case is treated):

\begin{definition}
A stationary process $Q$ on $[A^{\mathbb{Z}},\mathfrak{A}^{\mathbb{Z}}]$
will be called $\ast $-mixing if for each $0<\alpha <1$ there exists an $%
l\in \mathbb{N}$ such that 
\begin{equation}
\alpha Q(B)Q(C)\leq Q(B\cap C)\leq \alpha ^{-1}Q(B)Q(C)  \label{mix}
\end{equation}
whenever $B\in \mathfrak{A}^{\{...,-2,-1,0\}},C\in \mathfrak{A}%
^{\{l,l+1,...\}}$.
\end{definition}

Here $\mathfrak{A}^{T},T\subset \mathbb{Z}$, denotes the sub-$\sigma $-field
of $\mathfrak{A}^{\mathbb{Z}}$ concerning only times $t\in T$.

Observe that irreducible and aperiodic (i.e. weakly mixing) Markov chains
are automatically $\ast $-mixing, even with $\alpha =\alpha (l)$ tending to $%
1$ exponentially fast as $l\rightarrow \infty $. In the general situation,
even strong mixing ($\alpha $-mixing in the terminology of \cite{bradley})
does not imply $\ast $-mixing, because rare events may still deviate much
from independence.

We note that in the following we use the seemingly weaker condition, that (%
\ref{mix}) is fulfilled for \emph{some} $\alpha >0$ and \emph{some }$l$.
(The same was emphasized for most of the results in \cite{hiai-petz2}.) But,
in fact, in the stationary classical situation this is already equivalent to
full $\ast $-mixing, see \cite{bradley}, Theorem 4.1.

Let $\mathcal{P}_{\text{stat}}(A^{\mathbb{Z}})$, $\mathcal{P}_{\text{erg}%
}(A^{\mathbb{Z}})$ resp. $\mathcal{P}_{\ast }(A^{\mathbb{Z}})$ denote the
set of stationary, of ergodic resp. stationary $\ast $-mixing processes with
state space $A$.

We briefly introduce now the corresponding quantum set-up.

Consider a finite-dimensional $C^{\ast }$-algebra $\mathcal{A}$. The
classical case is covered choosing $\mathcal{A}$ is abelian.

It is well-known, that $\mathcal{A}$ can always be represented as a finite
direct sum of matrix algebras 
\begin{equation}
\mathcal{A}\cong \bigoplus_{i=1}^{m}\mathcal{M}_{k_{i}},  \label{reprmat}
\end{equation}
where $\mathcal{M}_{k}$ is the algebra of complex $k\times k$ matrices. The
abelian case is covered if all $k_{i}$ are $1$, meaning that $\mathcal{A}$
is simply the commutative algebra of complex functions over a finite set $%
A=\{1,2,...,m\}.$ A \emph{state }$\psi $ on $\mathcal{A}$ is a positive
functional on $\mathcal{A}$ with the property $\psi (\mathbf{1})=1$, where $%
\mathbf{1}$ is the unity. The set of all states on $\mathcal{A}$ is denoted
by $\mathcal{S}(\mathcal{A}).$ This is the set of probability measures on $A 
$ in the abelian case.

Any state $\psi $ on the finite-dimensional algebra $\mathcal{A}$ is
uniquely given by its \emph{density operator} $D_{\psi }\in \mathcal{A}$,
which is a positive trace-one operator fulfilling 
\begin{equation*}
\psi (X)=\text{Tr}_{\mathcal{A}}(D_{\psi }X)\text{ for each }X\in \mathcal{A}%
\text{.}
\end{equation*}
Here Tr$_{\mathcal{A}}$ denotes the canonical trace in $\mathcal{A}$ which
is nothing but the sum of the matrix traces in the above representation (\ref%
{reprmat}).

The quantum generalization of a stochastic process is usually constructed as
follows (and in correspondence to the definition of a process by its
compatible finite-dimensional distributions via Kolmogorov's extension
theorem): For each finite subset $T\subset \mathbb{Z}$ consider the $C^{\ast
}$-algebra $\mathcal{A}^{T}:=\bigotimes_{t\in T}\mathcal{A}$. Then for any $%
T\subset T^{\prime }\subset \mathbb{Z}$ there is a canonical embedding of $%
\mathcal{A}^{T}$ into $\mathcal{A}^{T^{\prime }}$ as a $C^{\ast }$%
-subalgebra. With respect to this identification consider the algebra 
\begin{equation*}
\widetilde{\mathcal{A}}:=\bigcup_{\substack{ T\subset \mathbb{Z}  \\ T\text{
finite}}}\mathcal{A}^{T}=\bigcup_{n\in \mathbb{N}}\mathcal{A}^{\{-n,...,n\}}.
\end{equation*}
$\widetilde{\mathcal{A}}$ is not norm-complete. We denote the completion by $%
\mathcal{A}^{\mathbb{Z}}$. It is a $C^{\ast }$-algebra and is called the 
\emph{quasilocal algebra} constructed from $\mathcal{A}$. Again, a state $%
\Psi $ on $\mathcal{A}^{\mathbb{Z}}$ is a positive functional on $\mathcal{A}%
^{\mathbb{Z}}$ with the property $\Psi (\mathbf{1})=1$. If $\mathcal{A}$ is
abelian, there is a one-to-one correspondence between states on $\mathcal{A}%
^{\mathbb{Z}}$ and stochastic processes with alphabet $A$: The restrictions $%
\Psi ^{(T)}:=\Psi \upharpoonright _{\mathcal{A}^{T}}$ of $\Psi $ to the
local algebras $\mathcal{A}^{T}$ correspond to the marginals $P^{(T)}$ of
the stochastic process $P$ on the cylinder $\sigma $-algebras $\mathfrak{A}%
^{T}$. This comes from the fact that any compatible family of local states $%
\Psi ^{(T)}$ has a unique extension to $\mathcal{A}^{\mathbb{Z}}$ just as
any compatible family of marginals can be extended to a stochastic process.

There is a canonically defined shift operator $\tau $ on $\mathcal{A}^{%
\mathbb{Z}}$ (mapping in particular $\mathcal{A}^{\{0\}}\subset \mathcal{A}^{%
\mathbb{Z}}$ onto $\mathcal{A}^{\{1\}}\subset \mathcal{A}^{\mathbb{Z}}$).
The set of \emph{stationary states} $\mathcal{S}_{\text{stat}}(\mathcal{A}^{%
\mathbb{Z}})$ is the subset of states in $\mathcal{S}(\mathcal{A}^{\mathbb{Z}%
})$ which are invariant with respect to $\tau $. This is a Choquet simplex,
the extremal points are called \emph{ergodic states }$\mathcal{S}_{\text{erg}%
}(\mathcal{A}^{\mathbb{Z}})$. The notions coincide with the classical ones
in the abelian case.

We complete the picture by defining a mixing property (cf. \cite{hiai-petz2}%
) as above:

\begin{definition}
\label{defn:ast-mixing} A stationary state $\Phi $ in $\mathcal{S}_{\text{%
stat}}(\mathcal{A}^{\mathbb{Z}})$ will be called $\ast $-mixing if for each $%
0<\alpha <1$ there exists an $l\in \mathbb{N}$ such that for each $k\in 
\mathbb{N}$%
\begin{eqnarray*}
&&\alpha \Phi ^{(\{-k,-k+1...,0\})}\otimes \Phi ^{(\{l,l+1,...,l+k\})} \\
&\leq &\Phi ^{(\{-k,-k+1...,0\}\cup \{l,l+1,...,l+k\})} \\
&\leq &\alpha ^{-1}\Phi ^{(\{-k,-k+1...,0\})}\otimes \Phi
^{(\{l,l+1,...,l+k\})}
\end{eqnarray*}
.
\end{definition}

We denote the set of stationary $\ast $-mixing states by $\mathcal{S}_{\ast
}(\mathcal{A}^{\mathbb{Z}})$.

Next we introduce the quantum version of the relative entropy rate. Let $%
\psi ,\varphi \in \mathcal{S}(\mathcal{A})$. The \emph{relative entropy }is
defined as 
\begin{equation*}
S(\psi ,\varphi ):=\left\{ 
\begin{array}{ll}
\text{Tr}_{\mathcal{A}}D_{\psi }(\log D_{\psi }-\log D_{\varphi }), & \text{%
if supp}(\psi )\leq \text{supp}(\varphi ) \\ 
\infty , & \text{otherwise.}%
\end{array}%
\right.
\end{equation*}%
Here supp$(D)$ is the smallest projection $p\in \mathcal{A}$ fulfilling $%
pDp=D$ (with $D\in \mathcal{A}$ self-adjoint).

Now, for $\Psi ,\Phi \in \mathcal{S}(\mathcal{A}^{\mathbb{Z}})$, we define
the relative entropy rate 
\begin{equation*}
s(\Psi ,\Phi ):=\lim_{n\rightarrow \infty }\frac{1}{n}S(\Psi ^{(n)},\Phi
^{(n)}),
\end{equation*}
whenever this limit exists in $\overline{\mathbb{R}}_{+}$ (we write for
short $\Psi ^{(n)}$ instead of $\Psi ^{(\{0,1,...,n-1\})}$ and $\mathcal{A}%
^{(n)}$ instead of $\mathcal{A}^{\{0,1,...,n-1\}}$).

Again, if $\Phi \in \mathcal{S}(\mathcal{A}^{\mathbb{Z}}),\Omega \subseteq 
\mathcal{S}(\mathcal{A}^{\mathbb{Z}})$ and $s(\Psi ,\Phi )$ exists for each $%
\Psi \in \Omega $ we write $s(\Omega ,\Phi )$ for $\inf_{\Psi \in \Omega
}s(\Psi ,\Phi )$.

\section{Equivalence of Sanov's theorem and Stein's lemma\label{chsteinsan}}

The maximally separating exponents for two states $\Psi ,\Phi $ on $\mathcal{%
A}^{\mathbb{Z}}$ are defined by 
\begin{equation*}
\beta _{\varepsilon ,n}(\Psi ,\Phi ):=\min \{\log \Phi (q):q\in \mathcal{A}%
^{(n)}\text{ projection, }\Psi (q)\geq 1-\varepsilon \},
\end{equation*}
for $\varepsilon \in (0,1)$. By $\bar{\beta}_{\varepsilon }(\Psi ,\Phi )$ we
denote limsup$_{n\rightarrow \infty }\frac{1}{n}\beta _{\varepsilon ,n}(\Psi
,\Phi )$, and if the limit exists in $-\overline{\mathbb{R}}_{+}:= - [0,
\infty]$, we denote it by $\beta _{\varepsilon }(\Psi ,\Phi )$.

\begin{definition}
We say that the pair $(\Psi ,\Phi )$ satisfies \emph{the HP-condition} if
the relative entropy rate $s(\Psi ,\Phi )$ exists and $\overline{\beta }%
_{\varepsilon }(\Psi ,\Phi )\leq -s(\Psi ,\Phi )$ for all $\varepsilon \in
(0,1)$.
\end{definition}

This condition was first proved to be fulfilled by Hiai and Petz in \cite%
{hiai-petz1} for the special case that $\Psi $ is \emph{completely ergodic}
and $\Phi $ is a \emph{stationary product state }(i.e. an iid state) and
later in \cite{hiai-petz2} for completely ergodic $\Psi $ and $\ast $-mixing
states $\Phi $.

\begin{definition}
We say that $\Phi \in \mathcal{S}(\mathcal{A}^{\mathbb{Z}})$ is a \emph{%
HP-state} if, for any ergodic state $\Psi \in \mathcal{S}_{\text{erg}}(%
\mathcal{A}^{\mathbb{Z}})$, the pair $(\Psi ,\Phi )$ satisfies the
HP-condition.
\end{definition}

As it turns out, the statement in Sanov's theorem is equivalent to the
HP-condition:

\begin{theorem}
\label{stein-sanov} Let $\Phi $ be a state on $\mathcal{A}^{\mathbb{Z}}$ and 
$\Theta \subseteq \mathcal{S}_{\text{erg}}(\mathcal{A}^{\mathbb{Z}})$. Then
following statements are equivalent:

\begin{itemize}
\item[1.] For each $\Psi \in \Theta $ the pair $(\Psi ,\Phi )$ satisfies the
HP-condition.

\item[2.] The quantity $s(\Psi ,\Phi )\leq +\infty $ exists for each $\Psi
\in \Theta $, and to each subset $\Omega \subseteq \Theta $ and any $\eta >0$
there exists a sequence $\{p_{n}\}_{n\in \mathbb{N}}$ of projections $%
p_{n}\in \mathcal{A}^{(n)}$ with 
\begin{equation}
\lim_{n\rightarrow \infty }\Psi ^{(n)}(p_{n})=1,\qquad \text{for all }\Psi
\in \Omega  \label{qqtyp}
\end{equation}%
such that if $s(\Omega ,\Phi )<\infty $ 
\begin{equation}
\underset{n\rightarrow \infty }{\text{\emph{limsup}}}\frac{1}{n}\log \Phi
^{(n)}(p_{n})\leq -s(\Omega ,\Phi )+\eta ,  \label{qqsep1}
\end{equation}%
otherwise 
\begin{equation}
\underset{n\rightarrow \infty }{\text{\emph{limsup}}}\frac{1}{n}\log \Phi
^{(n)}(p_{n})\leq -\frac{1}{\eta }.  \label{qqsep2}
\end{equation}%
Moreover, for each sequence of projections $\{\widetilde{p}_{n}\}$
fulfilling (\ref{qqtyp}) we have 
\begin{equation*}
\underset{n\rightarrow \infty }{\text{\emph{liminf}}}\frac{1}{n}\log \Phi
^{(n)}(\widetilde{p}_{n})\geq -s(\Omega ,\Phi )\text{.}
\end{equation*}%
Hence $-s(\Omega ,\Phi )$ is the lower limit of all achievable separation
exponents.
\end{itemize}
\end{theorem}

\begin{remark}

1. There are examples showing that in general one cannot choose $\eta =0$,
meaning that the exact value $-s(\Omega ,\Phi )$ is not necessarily
achievable.

2. If $\Phi $ is stationary and, moreover, $\ast $-mixing, statement 1. of
the Theorem is fulfilled with $\Theta =\mathcal{S}_{\text{erg}}(\mathcal{A}^{%
\mathbb{Z}})$. This will be seen in section \ref{starmix}.
\end{remark}

The implication $2.\Rightarrow 1.$ is trivial. The proof of the converse
implication is carried out in subsection \ref{prf}.

As an immediate consequence, we have the following assertion for the
classical case:

Let the maximally separating exponents for two processes $P,Q\in \mathcal{P}%
(A^{\mathbb{Z}})$ be defined by 
\begin{equation*}
\beta _{\varepsilon ,n}(P,Q):=\min \{\log Q^{(n)}(M):M\subseteq A^{n}\text{, 
}P^{(n)}(M)\geq 1-\varepsilon \},
\end{equation*}
for $\varepsilon \in (0,1)$. By $\bar{\beta}_{\varepsilon }(P,Q)$ we denote
limsup$_{n\rightarrow \infty }\frac{1}{n}\beta _{\varepsilon ,n}(P,Q)$, and
if the limit exists in $-\overline{\mathbb{R}}_{+}$, we denote it by $\beta
_{\varepsilon }(P,Q)$.

\begin{theorem}
\label{classsan}Let $Q$ $\in \mathcal{P}(A^{\mathbb{Z}}),\Theta \subseteq 
\mathcal{P}_{erg}(A^{\mathbb{Z}})$ and suppose that the relative entropy
rate $h(P,Q)$ exists for all $P\in \Theta $. Then following statements are
equivalent:

\begin{itemize}
\item[1.] $\bar{\beta}_{\varepsilon }(P,Q)\leq -h(P,Q)$ for all $P\in \Theta 
$ and all $\varepsilon \in (0,1)$.

\item[2.] For each set $\Omega \subseteq \Theta $ each $\eta >0$ there is a
sequence of subsets $\{M_{n}\}$, $M_{n}\subseteq A^{n}$, such that 
\begin{equation}
\lim_{n\rightarrow \infty }P(M_{n})=1,\qquad \text{for all }P\in \Omega ,
\label{classseparat}
\end{equation}
\end{itemize}

and 
\begin{equation*}
\underset{n\rightarrow \infty }{\text{\emph{limsup}}}\frac{1}{n}\log
Q^{(n)}(M_{n})\leq -h(\Omega ,Q)+\eta
\end{equation*}%
if $h(\Omega ,Q)<\infty $, otherwise if $h(\Omega ,Q)=\infty $ 
\begin{equation*}
\underset{n\rightarrow \infty }{\text{\emph{limsup}}}\frac{1}{n}\log
Q^{(n)}(M_{n})\leq -\frac{1}{\eta }.
\end{equation*}%
Moreover, for each sequence of subsets $\{\widetilde{M}_{n}\}$ fulfilling (%
\ref{classseparat}) we have 
\begin{equation*}
\underset{n\rightarrow \infty }{\text{\emph{liminf}}}\frac{1}{n}\log Q^{(n)}(%
\widetilde{M}_{n})\geq -h(\Omega ,Q).
\end{equation*}%
Hence $-h(\Omega ,Q)$ is the lower limit of all achievable separation
exponents.
\end{theorem}

\subsection{ A quantum relative AEP and achievability in Stein's lemma}

We start with a useful lemma which allows to translate some standard
techniques and estimates used in classical information theory into the
quantum setting.

\begin{lemma}
\label{mainest}Let $p,q$ be arbitrary projections and $\tau $ be a state on $%
\mathcal{A}$. Suppose that $u$ is a projection commuting with $D_{\tau }$.
Then we have 
\begin{equation}
\tau (qpq)\geq \tau (qpqu)\geq \tau (p)-2\left( \tau (\mathbf{1}-q)\right)
^{1/2}-\tau (\mathbf{1}-u) .  \label{psiqpqu}
\end{equation}
Let $c>0$. If $D_{\tau }u\leq cu$ then 
\begin{equation}
\text{\emph{Tr}}(pq)\geq \frac{1}{c}\left( \tau (p)-2\left( \tau (\mathbf{1}%
-q)\right) ^{1/2}-\tau (\mathbf{1}-u)\right) .  \label{traceppp}
\end{equation}
\end{lemma}

\begin{proof}
The first inequality in (\ref{psiqpqu}) is trivial. The second follows
applying the Cauchy-Schwarz inequality for the Hilbert-Schmidt inner
product: 
\begin{eqnarray*}
\tau (p) &=&\tau (pq)+\tau (p(\mathbf{1}-q)) \\
&\leq &|\tau (pq)|+(\text{Tr}(D_{\tau }(\mathbf{1}-q)))^{\frac{1}{2}} \\
&\leq &\tau (qpq)+|\tau ((\mathbf{1}-q)pq)|+\left( \tau (\mathbf{1}%
-q)\right) ^{1/2} \\
&\leq &\tau (qpq)+2\left( \tau (\mathbf{1}-q)\right) ^{1/2} \\
&=&\tau (qpqu)+\tau (qpq(\mathbf{1}-u))+2\left( \tau (\mathbf{1}-q)\right)
^{1/2} \\
&\leq &\tau (qpqu)+\tau (\mathbf{1}-u)+2\left( \tau (\mathbf{1}-q)\right)
^{1/2}.
\end{eqnarray*}
In the last inequality the assumption $[u, D_{\tau}]=0$ has been used.
Finally observe that $u\geq \frac{1}{c}D_{\tau }u$ and 
\begin{equation}
\text{Tr}(pq)=\text{Tr}(qpq)=\text{Tr}(qpqu)+\text{Tr}(qpq(\mathbf{1}%
-u))\geq \frac{1}{c}\text{Tr}(qpqD_{\tau }u).  \label{tracepq}
\end{equation}
Inequality (\ref{traceppp}) follows immediately inserting (\ref{psiqpqu})
into eqn. (\ref{tracepq}).$\qquad \qquad $
\end{proof}

For $0\leq s<\infty $, write $u_{\Phi ^{(n)}}^{\varepsilon }(s)$ for the
finite direct sum 
\begin{equation*}
u_{\Phi ^{(n)}}^{\varepsilon }(s):=\sum_{s-\varepsilon <s^{\prime
}<s+\varepsilon }\text{spec}_{e^{-ns^{\prime }}}(\Phi ^{(n)}),
\end{equation*}
where spec$_{\lambda }(\cdot )$ denotes the eigen-projection of its
argument's density operator (here $D_{\Phi ^{(n)}}$) corresponding to the
eigenvalue $\lambda $. We extend this definition to the case $s=\infty $ by
setting 
\begin{equation}
u_{\Phi ^{(n)}}^{\varepsilon }(\infty ):=\text{spec}_{0}(\Phi ^{(n)})+
\sum_{\varepsilon ^{-1}<s^{\prime }}\text{spec}_{e^{-ns^{\prime }}}(\Phi
^{(n)}).  \label{infrelent}
\end{equation}

Now we have

\begin{proposition}
\label{propaep}Let $\Psi $ be ergodic and let $\Phi $ be an arbitrary state
on $\mathcal{A}^{\mathbb{Z}}$. If the pair $(\Psi ,\Phi )$ fulfills the
HP-condition then:

\begin{itemize}
\item $\beta _{\varepsilon }(\Psi ,\Phi )=\lim_{n\rightarrow \infty }\frac{1%
}{n}\beta _{\varepsilon ,n}(\Psi ,\Phi )$ exists and we have 
\begin{equation*}
\beta _{\varepsilon }(\Psi ,\Phi )=-s(\Psi ,\Phi )
\end{equation*}
for each $\varepsilon \in (0,1)$.

\item Moreover, for all ${\varepsilon }>0$ it holds that 
\begin{equation*}
\lim_{n\rightarrow \infty }\Psi ^{(n)}(u_{\Phi ^{(n)}}^{\varepsilon }(s(\Psi
)+s(\Psi ,\Phi )))=1\text{ (relative AEP)}
\end{equation*}
and for each sequence $\{p_{n}\}$ of projections fulfilling 
\begin{equation}  \label{entropy-typical}
\Psi ^{(n)}(p_{n})\underset{n\rightarrow \infty }{\rightarrow }1\quad \text{
and }\quad \frac{1}{n}\log \text{\emph{Tr}}p_{n}\underset{n\rightarrow
\infty }{\rightarrow }s(\Psi )
\end{equation}
there is a sequence $\varepsilon _{n}\searrow 0$ such that with $%
u_{n}:=u_{\Phi ^{(n)}}^{\varepsilon _{n}}(s(\Psi )+s(\Psi ,\Phi ))$, the
relations 
\begin{equation*}
\Psi ^{(n)}(\text{\emph{supp}}(u_{n}p_{n}u_{n}))\underset{n\rightarrow
\infty }{\rightarrow }1
\end{equation*}
and 
\begin{equation*}
\frac{1}{n}\log \Phi ^{(n)}(\text{\emph{supp}}(u_{n}p_{n}u_{n}))\underset{%
n\rightarrow \infty }{\rightarrow }-s(\Psi ,\Phi )\text{ \ \ (max.
separating projection)}
\end{equation*}
are fulfilled.
\end{itemize}
\end{proposition}

\emph{Remark}. The quantum Shannon-McMillan theorem \cite{q-sm} guarantees
the existence of a sequence of projections $\{p_{n}\}$ with the properties
assumed in (\ref{entropy-typical}). We refer to such sequences as
entropy-typical w.r.t. $\Psi $. Roughly speaking, the above proposition
shows that one obtains a sequence of maximally seperating projections as an
'intersection' of $\Psi $-entropy-typical projections with appropriate
eigen-projections of the reference state $\Phi $.

\begin{proof}
1. First assume $s(\Psi ,\Phi )<+\infty $. By the monotonicity of the
relative entropy we may conclude that $S(\Psi ^{(n)},\Phi ^{(n)})<+\infty $
for each $n$. We have 
\begin{equation*}
\frac{1}{n}S(\Psi ^{(n)},\Phi ^{(n)})=-\frac{1}{n}S(\Psi ^{(n)})-\frac{1}{n}%
\text{Tr}D_{\Psi ^{(n)}}\log D_{\Phi ^{(n)}}\text{.}
\end{equation*}%
Let $\{\lambda _{i}\}_{i=1}^{r(\Phi ^{(n)})}$ be the set of non-zero
eigenvalues of $D_{\Phi ^{(n)}}$. We get 
\begin{equation*}
\frac{1}{n}S(\Psi ^{(n)},\Phi ^{(n)})=-\frac{1}{n}S(\Psi ^{(n)})-\frac{1}{n}%
\sum_{i}\log \lambda _{i}\text{Tr}D_{\Psi ^{(n)}}\text{spec}_{\lambda
_{i}}(\Phi ^{(n)}).
\end{equation*}%
Let ${\varepsilon }>0$ and 
\begin{equation*}
p_{n,{\varepsilon }}:=\sum_{\lambda _{_{i}}>e^{-n(s(\Psi )+s(\Psi ,\Phi
)-\varepsilon )}}\text{spec}_{\lambda _{i}}(\Phi ^{(n)}).
\end{equation*}%
We claim that 
\begin{equation*}
\lim_{n\rightarrow \infty }\Psi (p_{n,{\varepsilon }})=0\qquad \text{for all 
}{\varepsilon }>0.
\end{equation*}%
In fact, suppose on the contrary that for some $\varepsilon >0$ we have 
\begin{equation*}
\underset{n\rightarrow \infty }{\text{limsup}}\Psi (p_{n,{\varepsilon }})>0.
\end{equation*}%
We conclude the existence of some $\gamma >0$ and some subsequence $%
\{n_{j}\} $ with 
\begin{equation*}
\Psi ^{(n_{j})}(p_{n_{j},{\varepsilon }})>\gamma >0.
\end{equation*}%
Fix some $\alpha \in (0,1),\delta >0$. Let $p_{n_{j}}:=p_{n_{j},\varepsilon
} $, $q_{n_{j}}:=\arg \min \beta _{\alpha ,n_{j}}(\Psi ,\Phi )$ and $%
u_{n_{j}}:=u_{\Psi ^{(n_{j})}}^{\delta }(s(\Psi ))$. Then $D_{\Psi
^{(n_{j})}}u_{n_{j}}\leq cu_{n_{j}}$ for $c=e^{-n_{j}(s(\Psi )-\delta )}$
and by Lemma \ref{mainest} and the quantum Shannon-McMillan theorem we
arrive at 
\begin{equation*}
\Psi ^{(n_{j})}(q_{n_{j}}p_{n_{j}}q_{n_{j}}u_{n_{j}})\geq \gamma -2\sqrt{%
\alpha }-\delta >0,
\end{equation*}%
and 
\begin{equation}
\text{Tr}(q_{n_{j}}p_{n_{j}}q_{n_{j}})\geq e^{n_{j}(s(\Psi )-\delta
)}(\gamma -2\sqrt{\alpha }-\delta ),  \label{firstx}
\end{equation}%
if $j$ is large enough and if $2\sqrt{\alpha }+\delta <\gamma $. Now,
observe that $D_{\Phi ^{(n_{j})}}$ and $p_{n_{j}}$ commute and that
consequently we have $D_{\Phi ^{(n_{j})}}\geq e^{-n_{j}(s(\Psi )+s(\Psi
,\Phi )-\varepsilon )}p_{n_{j}}$ by definition of $p_{n_{j}}$. Thus we
obtain 
\begin{equation}
q_{n_{j}}D_{\Phi ^{(n_{j})}}q_{n_{j}}\geq e^{-n_{j}(s(\Psi )+s(\Psi ,\Phi
)-\varepsilon )}q_{n_{j}}p_{n_{j}}q_{n_{j}}.  \label{secondx}
\end{equation}%
After applying trace to both sides of this inequality, taking logarithms,
dividing by $n_{j}$, taking limit superior and using (\ref{firstx}) we are
led to 
\begin{equation*}
\bar{\beta}_{\alpha }(\Psi ,\Phi )\geq -s(\Psi ,\Phi )+{\varepsilon }-\delta
>-s(\Psi ,\Phi ),
\end{equation*}%
which contradicts the assumed HP-condition provided that $\delta <{%
\varepsilon }$. \newline
In the case $s(\Psi ,\Phi )=\infty $ everything can be done in the same way,
we just have to substitute the definition of $p_{n,{\varepsilon }}$ by 
\begin{equation*}
p_{n,{\varepsilon }}=\sum_{\lambda _{_{i}}>e^{-n/\varepsilon }}\text{spec}%
_{\lambda _{i}}(\Phi ^{(n)})
\end{equation*}%
\ and obtain $\bar{\beta}_{\alpha }\geq -\frac{1}{{\varepsilon }}+s(\Psi
)-\delta $, again in contradiction to $\bar{\beta}_{\alpha }(\Psi ,\Phi
)=-\infty $, hence again the projectors $p_{n,{\varepsilon }}$ have
asymptotically vanishing expectation with respect to $\Psi $ for each
positive $\varepsilon $.

2. Let first $s(\Psi ,\Phi )<\infty $. We have $\frac{1}{n}S(\Psi
^{(n)},\Phi ^{(n)})\rightarrow s(\Psi ,\Phi )$ as $n\rightarrow \infty $ by
assumption, hence 
\begin{equation*}
-\frac{1}{n}\text{Tr}D_{\Psi ^{(n)}}\log D_{\Phi ^{(n)}}\underset{%
n\rightarrow \infty }{\longrightarrow }s(\Psi )+s(\Psi ,\Phi )
\end{equation*}
and the mixed term $-\frac{1}{n}$Tr$D_{\Psi ^{(n)}}\log D_{\Phi ^{(n)}}$ is
the expectation value of the random variable $-\frac{1}{n}\log \lambda _{i}$
with respect to the probability measure given by 
\begin{equation*}
\{\text{Tr}D_{\Psi ^{(n)}}\text{spec}_{\lambda _{i}}(\Phi ^{(n)})\},
\end{equation*}
where again $\{\lambda _{i}\}$ runs through the non-zero eigenvalues of $%
\Phi ^{(n)}$. On the other hand, we have shown in 1. that the \emph{lower
bounded }random variable $-\frac{1}{n}\log \lambda _{i}\geq 0$ is bounded
asymptotically in probability by the quantity $s(\Psi )+s(\Psi ,\Phi )$,
being its asymptotic expectation value at the same time, i.e $%
\lim_{n\rightarrow \infty }\Psi (p_{n,{\varepsilon }})=0$ for all ${%
\varepsilon }>0$. From this it easily follows that 
\begin{equation*}
\lim_{n\rightarrow \infty }\Psi (t_{n,\delta })=0\quad \text{for all }\quad
\delta >0,
\end{equation*}
where 
\begin{equation*}
t_{n,\delta }:=\sum_{\lambda _{i}<e^{-n(s(\Psi )+s(\Psi ,\Phi )+\delta )}}%
\text{spec}_{\lambda _{i}}(\Phi ^{(n)}).
\end{equation*}
This is the assertion 
\begin{equation}
\lim_{n\rightarrow \infty }\Psi ^{(n)}(u_{\Phi ^{(n)}}^{\varepsilon }(s(\Psi
)+s(\Psi ,\Phi )))=1  \label{typicspec}
\end{equation}
for all $\varepsilon >0$. In the case $s(\Psi ,\Phi )=\infty $ the relative
AEP follows immediately from 1.

3. First assume $s(\Psi ,\Phi )<\infty $. Fix some $\varepsilon $ and some $%
\alpha \in (0,1)$. Let $\{q_{n}\}$ be any sequence of projections fulfilling 
$\Psi ^{(n)}(q_{n})\geq 1-\alpha $ for $n$ large enough. Let $p_{n}:=u_{\Phi
^{(n)}}^{\varepsilon }(s(\Psi )+s(\Psi ,\Phi ))$. We proved that $\Psi
^{(n)}(p_{n})\underset{n\rightarrow \infty }{\rightarrow }1$. Now as in 1.
we may conclude 
\begin{equation*}
\Phi ^{(n)}(q_{n})\geq e^{-n(s(\Psi )+s(\Psi ,\Phi )+\varepsilon )}\text{Tr}%
p_{n}q_{n}\text{.}
\end{equation*}%
Using the quantum Shannon-McMillan theorem and again Lemma \ref{mainest},
this time applied to the density operator of $\Psi ^{(n)}$ and with $%
u:=\sum_{\lambda <e^{-n(s(\Psi )-\delta )}}$spec$_{\lambda }(\Psi ^{(n)})$,
for arbitrary $\delta >0$ and $n$ large enough we have Tr$p_{n}q_{n}\geq
e^{n(s(\Psi )-\delta )}a$ for some $0<a<1$ independent of $n$. Hence we get
for any $\varepsilon ,\delta >0$%
\begin{equation*}
\frac{1}{n}\log \Phi ^{(n)}(q_{n})>-s(\Psi ,\Phi )-\varepsilon -\delta
\end{equation*}%
for $n$ large enough. Therefore, the quantity $\beta _{\alpha }(\Psi ,\Phi )$
exists for any $\alpha \in (0,1)$ and coincides with $-s(\Psi ,\Phi )$,
where we used again the HP-condition. In the case $s(\Psi ,\Phi )=\infty $
this assertion is a trivial consequence of $\bar{\beta}_{\alpha }(\Psi ,\Phi
)=-\infty $.

4. Let $\{p_{n}\}$ be a sequence of projections $p_{n}\in \mathcal{A}^{(n)}$
with $\lim_{n\rightarrow \infty }\Psi ^{(n)}(p_{n})=1$ and $%
\lim_{n\rightarrow \infty }\frac{1}{n}\log \text{Tr}(p_{n})=s(\Psi )$. Fix
some $\varepsilon >0$. Let us write $u_{n}$ instead of $u_{\Phi
^{(n)}}^{\varepsilon }(s(\Psi )+s(\Psi ,\Phi ))$ for short. From (\ref%
{typicspec}) and Lemma \ref{mainest} we infer that $\Psi
^{(n)}(u_{n}p_{n}u_{n})\underset{n\rightarrow \infty }{\rightarrow }1$. 
Now $u_{n}p_{n}u_{n}$ is a positive operator being upper bounded by its
support projection supp$(u_{n}p_{n}u_{n})$ which proves $\Psi ^{(n)}($supp$%
(u_{n}p_{n}u_{n}))\underset{n\rightarrow \infty }{\rightarrow }1$. From this
we easily conclude that we may even substitute the $\varepsilon $ in the
definition of $u_{n}=u_{\Phi ^{(n)}}^{\varepsilon }(s(\Psi )+s(\Psi ,\Phi ))$
by a suitable sequence $\varepsilon _{n}\rightarrow 0$ and still have $\Psi
^{(n)}($supp$(u_{n}p_{n}u_{n}))\underset{n\rightarrow \infty }{\rightarrow }%
1 $.\newline
On the other hand, we have $\text{supp}(u_{n}p_{n}u_{n})\leq u_{n}$ as well
as $\text{Tr}(\text{supp}(u_{n}p_{n}u_{n}))\leq \text{Tr}(p_{n})$. Hence we
get in the case $s(\Psi ,\Phi )<\infty $ 
\begin{eqnarray*}
&&\frac{1}{n}\log \Phi ^{(n)}(\text{supp}(u_{n}p_{n}u_{n})) \\
&\leq &\frac{1}{n}(-n(s(\Psi )+s(\Psi ,\Phi )-\varepsilon _{n})+\log \text{Tr%
}(p_{n}))\underset{n\rightarrow \infty }{\longrightarrow }-s(\Psi ,\Phi )%
\text{,}
\end{eqnarray*}
resp. for $s(\Psi ,\Phi )=\infty $ 
\begin{eqnarray*}
&&\frac{1}{n}\log \Phi ^{(n)}(\text{supp}(u_{n}p_{n}u_{n})) \\
&\leq &\frac{1}{n}(-n/\varepsilon _{n}+\log \text{Tr}(p_{n}))\underset{%
n\rightarrow \infty }{\longrightarrow }-s(\Psi ,\Phi )=-\infty \text{,}
\end{eqnarray*}

This together with the fact we proved that no sequence of $\Psi $ -typical
projections has a better lower limit of the separation rate than $-s(\Psi
,\Phi )$ shows now that 
\begin{equation*}
\frac{1}{n}\log \Phi ^{(n)}(\text{supp}(u_{n}p_{n}u_{n}))\underset{%
n\rightarrow \infty }{\longrightarrow }-s(\Psi ,\Phi )\text{.}
\end{equation*}
We proved all assertions of the proposition.$\qquad \qquad $ 
\end{proof}

\subsection{Stein's Lemma implies Sanov's Theorem\label{prf}}

With the preliminaries given in the last subsection, it is now easy to
complete the proof of Theorem \ref{stein-sanov}:

\begin{proof}
Let 
\begin{equation*}
s_{\text{min}}:=\inf_{\Psi \in \Omega }s(\Psi )\quad \text{and}\quad s_{%
\text{max}}:=\sup_{\Psi \in \Omega }s(\Psi ),
\end{equation*}
where $s(\Psi )$ denotes the von Neumann entropy rate of the ergodic state $%
\Psi $ $\in \Omega $. Choose $s_{1},\ldots ,s_{m}$ satisfying (for $\eta
:=m^{-1}(s_{\text{max}}-s_{\text{min}})$) 
\begin{equation*}
s_{\text{min}}=s_{1}<s_{2}<\ldots s_{m-1}<s_{m}=s_{\text{max}}\quad \text{and%
}\quad s_{i}-s_{i-1}=\eta ,i\in \{2,\ldots ,m\}.
\end{equation*}
Define $s_{m+1}=s_{m}+\eta $.

Let first $s(\Omega ,\Phi )<+\infty $. For $i\in \{1,\ldots ,m\}$ we
consider the collection of disjoint intervals 
\begin{equation*}
I_{i}:=\left( s_{i}+s(\Omega ,\Phi )-\frac{\eta }{2},s_{i}+s(\Omega ,\Phi )+%
\frac{\eta }{2}\right]
\end{equation*}
and 
\begin{eqnarray*}
I_{m+1}:=(s_{\text{max}}+s(\Omega, \Phi) + \frac{\eta}{2}, \infty).
\end{eqnarray*}
Moreover we define the following projections 
\begin{equation*}
u_{n,i}:=\sum_{-\frac{1}{n}\log \lambda \in I_{i}}\text{spec}_{\lambda
}(\Phi ^{(n)}),
\end{equation*}
and 
\begin{equation*}
u_{n,m+1}:=\sum_{-\frac{1}{n}\log \lambda > s_{\text{max}}+s(\Omega ,\Phi
)+\eta /2} \text{spec}_{\lambda }(\Phi ^{(n)}),
\end{equation*}
where the summations extend over the eigenvalues of $D_{\Phi ^{(n)}}$.
Additionally we consider universally typical projections $p_{n,i},i\in
\{1,\ldots ,m\}$ (according to the Kaltchenko-Yang universality result \cite%
{ky})\ to the levels $s_{i}+\eta $ (i.e. 
\begin{equation}
\lim_{n\rightarrow \infty }\frac{1}{n}\log \text{Tr}(p_{n,i})=s_{i}+\eta
\label{lalala1}
\end{equation}
and 
\begin{equation}
\lim_{n\rightarrow \infty }\Psi^{(n)} (p_{n,i})=1  \label{lalala2}
\end{equation}
for each ergodic state with $s(\Psi )<s_{i}+\eta $). In addition, set $%
p_{n,m+1}:=p_{n,m}$. We may choose the sequence of these projections to be
ascending, i.e. 
\begin{equation}
p_{n,i}\leq p_{n,i+1},  \label{ascend}
\end{equation}
since otherwise we may define 
\begin{equation*}
\widehat{p}_{n,i}:=\bigvee_{j=1}^{i}p_{n,j}.
\end{equation*}
The $\widehat{p}_{n,i}$ fulfil (\ref{lalala1}) and (\ref{lalala2}) as well,
so we may work with these instead of $p_{n,i}$.

Set $r_{n,i}=$supp$(u_{n,i}p_{n,i}u_{n,i})$ for $i=1,2,...,m+1$ and define $%
p_{n}$ by 
\begin{equation*}
p_{n}:=\sum_{i=1}^{m+1}r_{n,i}.
\end{equation*}
(Observe that the $r_{n,i}$ are mutually orthogonal.)

For $\Psi \in \Omega $ let $i_{0}\in \{1, \dots, m+1\}$ be the index
fulfilling $s(\Psi )+s(\Psi ,\Phi )\in I_{i_{0}}$. This means that $s(\Psi
)\leq s_{i_{0}}+\eta /2<s_{i_{0}}+\eta $. Consequently, by (\ref{lalala2})
we obtain $\lim_{n\rightarrow \infty }\Psi^{(n)} (p_{n,i_{0}})=1$.

Further, by the relative AEP (Proposition \ref{propaep}) $\lim_{n\rightarrow
\infty }\Psi^{(n)} (u_{n,i_{0}}+u_{n,i_{0}+1})=1$, for $i_0\in \{1, \dots,
m\}$, and $\lim_{n\rightarrow \infty }\Psi^{(n)} (u_{n,m+1})=1$ are
satisfied. We add the projection $u_{n,i_{0}+1}$ for $i_0 \in \{1, \dots,
m\} $ in order to cover the case where the mixed term is equal to the right
end point of $I_{i_0}$. We conclude from (\ref{ascend}) and Lemma \ref%
{mainest} that $\Psi^{(n)} (r_{n,i_{0}}+r_{n,i_{0}+1})\rightarrow 1$, for $%
i_0 \in \{1, \dots, m\}$, and $\Psi^{(n)} (r_{n,m+1})\rightarrow 1$.
Therefore 
\begin{equation*}
\lim_{n\rightarrow \infty }\Psi^{(n)} (p_{n})=1.
\end{equation*}

On the other hand we have for $n$ sufficiently large by (\ref{lalala1}) and
by definition of $\eta $ 
\begin{eqnarray*}
\Phi (p_{n}) &=&\sum_{i=1}^{m+1}\Phi (r_{n,i})\leq \sum_{i=1}^{m+1}\text{Tr}%
(p_{n,i})e^{-n(s_{i}+s(\Omega ,\Phi )-\eta /2)} \\
&\leq &\sum_{i=1}^{m+1}e^{n(s_{i}+2\eta )}e^{-n(s_{i}+s(\Omega ,\Phi )-\eta
/2)} \\
&=&e^{-n(s(\Omega ,\Phi )-\frac{5}{2}\eta -\frac{\log (m+1)}{n})} \\
&=&e^{-n(s(\Omega ,\Phi )-\frac{5}{2}m^{-1}(s_{\text{max}}-s_{\text{min}})-%
\frac{\log (m+1)}{n})}.
\end{eqnarray*}
So, by choosing $m$ sufficiently large, we get statement (\ref{qqsep1}).

The case $s(\Omega ,\Phi )=+\infty $ easily follows by setting $p_n:=u_{\Phi
^{(n)}}^{\eta }(\infty )$, see (\ref{infrelent}). By Proposition \ref%
{propaep} the projection $p_{n}$ is asymptotically typical for all $\Psi \in
\Omega $, and we have $\Phi ^{(n)}(p_{n})\leq $Tr$(\mathbf{1}_{\mathcal{A}%
^{(n)}})e^{-n\eta ^{-1}}=e^{-n(\eta ^{-1}-\log \text{Tr}(\mathbf{1}_{%
\mathcal{A}}))}$, so again we get statement (\ref{qqsep2}).

Finally, the fact that a better separation exponent than $-s(\Omega ,\Phi )$
is not achievable immediately follows from Proposition \ref{propaep}. $%
\qquad \qquad $
\end{proof}

\section{$\ast $-mixing implies the HP-condition\label{starmix}}

We start with a proposition extending the result in \cite{hiai-petz2} to the
case of only ergodic (instead of completely ergodic) $\Psi $.

\begin{theorem}
Let $\Phi \in \mathcal{S}_{\ast }(\mathcal{A}^{\mathbb{Z}})$. Then $\Phi $
is an HP-state.
\end{theorem}

Recall that $\mathcal{S}_{\ast }(\mathcal{A}^{\mathbb{Z}})$ denotes the set
of stationary $\ast $-mixing states, see Definition \ref{defn:ast-mixing}. 
\newline
\newline

\begin{proof}
1. Let $\Psi \in \mathcal{S}_{\text{erg}}(\mathcal{A}^{\mathbb{Z}})$. The
relative entropy rate $s(\Psi ,\Phi )\leq +\infty $ exists in view of \cite%
{hiai-petz2}, Theorem 2.1, in connection with Remark 4.2, \emph{ibid.} (even
if only $\Psi \in \mathcal{S}_{\text{stat}}(\mathcal{A}^{\mathbb{Z}})$\ is
assumed).

2. Fix some $l\in \mathbb{N}$, another integer $m$ (which in the sequel has
to be chosen large enough) and represent the quasilocal $C^{\ast }$-algebra $%
\mathcal{A}^{\mathbb{Z}}$ as $C^{\ast }$-algebra $(\mathcal{A}^{\otimes
l}\otimes \mathcal{A}^{\mathbb{\otimes }m})^{\mathbb{Z}}$, i.e. partition
the integers into blocks of length $l+m$, where each block consists of a
starting part of length $l$ and the remaining part of length $m$. Clearly,
the entropy rate $s_{(l+m)}(\Psi ,\Phi )$ with respect to this new
partitioning exists, and we have $s_{(l+m)}(\Psi ,\Phi )=(l+m)s(\Psi ,\Phi )$%
. With respect to the canonical shift operator $\tau _{l,m}:=\tau ^{l+m}$
acting in $(\mathcal{A}^{\otimes l}\otimes \mathcal{A}^{\mathbb{\otimes }%
m})^{\mathbb{Z}}$, the state $\Psi $ is still stationary, but may fail to be
ergodic. Anyway, it has a finite ergodic decomposition 
\begin{equation*}
\Psi =\frac{1}{l+m}\sum_{r=0}^{l+m-1}\Psi _{(r,l+m)},
\end{equation*}
where some of the ergodic components may coincide, and all $\Psi _{(r,l+m)}$
have the same entropy rate $s_{(l+m)}(\Psi _{(r,l+m)})\equiv s_{(l+m)}(\Psi
)=(l+m)s(\Psi )$, \cite{q-sm}. The ergodic components also have the same
relative entropy rate $s_{(l+m)}(\Psi _{(r,l+m)},\Phi ):= s_{(l+m)}(\Psi
,\Phi )=(l+m)s(\Psi ,\Phi )$. Observe that this was shown in \cite{igrai}
for the case of stationary product states $\Phi $. However as the proof only
makes use of the existence of the relative entropy rates $s_{(l+m)}(\Psi
_{(r,l+m)},\Phi )$, which is guaranteed in our situation (see 1.), the
monotonicity of the relative entropy and the affinity of the relative
entropy rate with respect to its first argument the relation extends to $%
\ast $-mixing reference states.

Next, denote by $\mathbf{I}$ the trivial subalgebra of $\mathcal{A}$
generated by the unit element $\mathbf{1}_{\mathcal{A}}$\ and\ consider the $%
C^{\ast }$-subalgebra $\mathbf{I}^{\otimes l}\otimes \mathcal{A}^{\mathbb{%
\otimes }m}$\ of $\mathcal{A}^{\otimes l}\otimes \mathcal{A}^{\mathbb{%
\otimes }m}$. Then $(\mathbf{I}^{\otimes l}\otimes \mathcal{A}^{\mathbb{%
\otimes }m})^{\mathbb{Z}}$ is a $C^{\ast }$-subalgebra of the quasi-local
algebra $(\mathcal{A}^{\otimes l}\otimes \mathcal{A}^{\mathbb{\otimes }m})^{%
\mathbb{Z}}$, and by \cite{bratrob}, Theorem 4.3.17, the restrictions $%
\widehat{\Psi }_{(r,l,m)}$ of the ergodic components $\Psi _{(r,l+m)}$ to $(%
\mathbf{I}^{\otimes l}\otimes \mathcal{A}^{\mathbb{\otimes }m})^{\mathbb{Z}}$
are ergodic, too. They are the ergodic components of $\widehat{\Psi }%
_{(l,m)}:=\Psi \upharpoonright _{(\mathbf{I}^{\otimes l}\otimes \mathcal{A}^{%
\mathbb{\otimes }m})^{\mathbb{Z}}}$.

We introduce the ($\tau _{l,m}$-)stationary product state $\widetilde{\Phi }%
_{(l+m)}$ on $(\mathcal{A}^{\otimes l}\otimes \mathcal{A}^{\mathbb{\otimes }%
m})^{\mathbb{Z}}$ which is uniquely defined by its one-site marginals $\Phi
\upharpoonright _{\mathcal{A}^{\otimes l}\otimes \mathcal{A}^{\mathbb{%
\otimes }m}}$, and consider its restriction $\widehat{\Phi }_{(l,m)}$\ to
the $C^{\ast }$-algebra$\ (\mathbf{I}^{\otimes l}\otimes \mathcal{A}^{%
\mathbb{\otimes }m})^{\mathbb{Z}},$ which is a ($\tau _{l,m}$-)stationary
product state, too.

3. In the following we have to take into account whether $s(\Psi ,\Phi )$ is
finite or infinite. Let us first treat the case $s(\Psi ,\Phi )<+\infty $.

Define for two states $\psi ,\varphi $ on a $C^{\ast }$-algebra $\mathcal{A}$%
, 
\begin{equation*}
S_{\text{co}}(\psi ,\varphi ):=\sup \left\{ \sum_{i}\psi (q_{i})\frac{\log
\psi (q_{i})}{\log \varphi (q_{i})}:q_{i}\text{ projections with }\sum
q_{i}=I\right\}
\end{equation*}
(cf. \cite{hiai-petz1}).

Consider the relative entropies $S(\widehat{\Psi }_{(r,l,m)}^{(l+m)},%
\widehat{\Phi }_{(l,m)}^{(l+m)})$, then we get by the subadditivity of
entropy and by the fact that the rates of the quantities $S$ and $S_{\text{co%
}}$ coincide (Hiai and Petz \cite{hiai-petz1}) 
\begin{eqnarray}
&&S(\widehat{\Psi }_{(r,l,m)}^{(l+m)},\widehat{\Phi }_{(l,m)}^{(l+m)})
\label{rrrrrr} \\
&\leq &s_{(l+m)}(\widehat{\Psi }_{(r,l,m)},\widehat{\Phi }%
_{(l,m)})=\lim_{k\rightarrow \infty }\frac{1}{k}S_{\text{co}}(\widehat{\Psi }%
_{(r,l,m)}^{(k(l+m))},\widehat{\Phi }_{(l,m)}^{(k(l+m))}).  \notag
\end{eqnarray}%
From the definition of $S_{\text{co}}$ and $\ast $-mixing we obtain now 
\begin{eqnarray*}
&&S(\widehat{\Psi }_{(r,l,m)}^{(l+m)},\widehat{\Phi }_{(l,m)}^{(l+m)}) \\
&\leq &\lim_{k\rightarrow \infty }\frac{1}{k}\left( S_{\text{co}}(\widehat{%
\Psi }_{(r,l,m)}^{(k(l+m))},\Phi \upharpoonright _{(\mathbf{I}^{\otimes
l}\otimes \mathcal{A}^{\mathbb{\otimes }m})^{\otimes k}})-k\log \alpha
\right) .
\end{eqnarray*}%
This is the same technique as used by Hiai and Petz in \cite{hiai-petz2}.
Again from the definition of $S_{\text{co}}$ we get 
\begin{eqnarray*}
&&S(\widehat{\Psi }_{(r,l,m)}^{(l+m)},\widehat{\Phi }_{(l,m)}^{(l+m)}) \\
&\leq &-\log \alpha +\lim_{k\rightarrow \infty }\frac{1}{k}S_{\text{co}%
}(\Psi _{(r,l+m)}^{(k(l+m))},\Phi ^{(k(l+m))}),
\end{eqnarray*}%
and from the relation $S_{\text{co}}\leq S$ (see \cite{hiai-petz1}), which
is a consequence of the monotonicity of the relative entropy, we arrive at 
\begin{eqnarray}
&&S(\widehat{\Psi }_{(r,l,m)}^{(l+m)},\widehat{\Phi }_{(l,m)}^{(l+m)})
\label{upest} \\
&\leq &-\log \alpha +\lim_{k\rightarrow \infty }\frac{1}{k}S(\Psi
_{(r,l+m)}^{(k(l+m))},\Phi ^{(k(l+m))})  \notag \\
&\leq &-\log \alpha +s_{(l+m)}(\Psi _{(r,l+m)},\Phi )=-\log \alpha
+(l+m)s(\Psi ,\Phi ).  \notag
\end{eqnarray}

This upper bound may be utilized to derive an essential lower bound. For an
arbitrarily chosen $\eta >0$, define 
\begin{equation*}
A_{l,m,\eta }:=\{r:0\leq r<l+m,\frac{1}{l+m}S(\widehat{\Psi }%
_{(r,l,m)}^{(l+m)},\widehat\Phi^{(m+l)}_{(m,l)})<s(\Psi ,\Phi )-\eta \}.
\end{equation*}
The convexity of the relative entropy in its first argument together with (%
\ref{upest})\ yields 
\begin{eqnarray}
\frac{1}{l+m}S(\Psi ^{(m)},\Phi ^{(m)}) &=&\frac{1}{l+m}S(\widehat\Psi
^{(m+l)}_{(m,l)},\widehat\Phi ^{(m+l)}_{(l,m)})  \label{lowest1} \\
&\leq &\frac{1}{(l+m)^{2}}\sum_{r}S(\widehat{\Psi }_{(r,l,m)}^{(l+m)},%
\widehat\Phi^{(m+l)}_{(m,l)})  \notag \\
&<&\frac{\#A_{l,m,\eta }}{l+m}(s(\Psi ,\Phi )-\eta )+\frac{\#A_{l,m,\eta
}^{c}}{l+m}\left( -\frac{\log \alpha }{l+m}+s(\Psi ,\Phi )\right) .  \notag
\end{eqnarray}
Fixing $l$ and letting $m\rightarrow \infty $, the expression $\frac{1}{l+m}%
S(\Psi ^{(m)},\Phi ^{(m)})$ tends to $s(\Psi ,\Phi )$. This immediately
leads to the conclusion that 
\begin{equation}
\frac{\#A_{l,m,\eta }}{l+m}\underset{m\rightarrow \infty }{\longrightarrow }0%
\text{ \ \ for each }l,\eta \text{.}  \label{denslem}
\end{equation}

For each $r\in A_{l,m,\eta }^{c}$, the relative entropy rate fulfills 
\begin{equation}
\frac{1}{l+m}s_{(l+m)}(\widehat{\Psi }_{(r,l,m)},\widehat{\Phi }%
_{(l,m)})\geq s(\Psi ,\Phi )-\eta ,  \label{relentlow}
\end{equation}
too, since $\widehat{\Phi }_{(l,m)}$ is a $\tau _{l+m}$-stationary product
state.

4. Hence we are in the situation treated in \cite{igrai}. The main assertion
of \cite{igrai} is the quantum Stein's lemma saying that for any given $%
\varepsilon >0$ it is possible to construct projections $p_{r,n,\varepsilon
}\in (\mathbf{I}^{\otimes l}\otimes \mathcal{A}^{\mathbb{\otimes }%
m})^{\otimes n}$ which are $\varepsilon $-typical with respect to $\widehat{%
\Psi }_{(r,l,m)}$ (i.e. $\widehat{\Psi }_{(r,l,m)}^{(n(l+m))}(p_{r,n,%
\varepsilon })\geq 1-\varepsilon $ for large $n$) and maximally separating: $%
\widehat{\Phi }_{(l,m)}^{((l+m)n)}(p_{r,n,\varepsilon })\leq e^{-n(s_{(l+m)}(%
\widehat{\Psi }_{(r,l,m)},\widehat{\Phi }_{(l,m)})-\varepsilon )}$ for large 
$n$. Moreover, the quantum relative AEP (Theorem 2 in \cite{igrai}) ensures,
in particular, that if $n$ is sufficiently large we have 
\begin{eqnarray*}
\text{Tr}p_{r,n,\varepsilon } &\leq &e^{n(s_{(l+m)}(\widehat{\Psi }%
_{(r,l,m)})+\varepsilon )}\text{ \ \ \ \ \ \ \ \ \ \ \ \ and} \\
\widehat{\Phi }_{(l,m)}^{((l+m)n)}(p) &\leq &e^{-n(s_{(l+m)}(\widehat{\Psi }%
_{(r,l,m)})+s_{(l+m)}(\widehat{\Psi }_{(r,l,m)},\widehat{\Phi }%
_{(l,m)})-\varepsilon )}
\end{eqnarray*}%
for each minimal projection $p\leq p_{r,n,\varepsilon }$.

For our purpose we need a bit more information about the construction of
these maximally separating projections. In the course of the proof of
Theorem 2 in \cite{igrai} the projections $p_{r,n,\varepsilon }$ are
constructed in the following way:

a) A super-block length $L$ is chosen, where the only requirement about $L$
is that it is large enough to ensure some appropriate entropy approximation
(any larger $L$ will do, too).

b) The projections $p_{r,nL,\varepsilon }$ are constructed as certain
sub-projections of the projection 
\begin{equation}
p^{(nL)}:=\sum_{\substack{ -\frac{1}{nL}\log \lambda \geq s_{(l+m)}(\widehat{%
\Psi }_{(r,l,m)})  \\ +s_{(l+m)}(\widehat{\Psi }_{(r,l,m)},\widehat{\Phi }%
_{(l,m)})-\varepsilon }}\text{spec}_{\lambda }((\widehat{\Phi }%
_{(l,m)}^{(l+m)})^{\otimes nL}).  \label{specproject}
\end{equation}

c) The remaining projections $p_{r,nL+k,\varepsilon }$ for $1\leq k<L$ are
constructed as $p_{r,nL,\varepsilon }\otimes (I_{\mathcal{A}^{\otimes
(l+m)}})^{\otimes k}$.

For given $l$ and $m$ we may choose one and the same super-block length $L$
for the different $r\in \{0,1,...,l+m-1\}$ and define our separating
projections first for the multiples of $L(l+m)$ by$\ $%
\begin{equation*}
q_{nL(l+m),\varepsilon }:=\bigvee_{r\in A_{l,m,\varepsilon
}^{c}}p_{r,nL,\varepsilon }.
\end{equation*}%
In view of (\ref{specproject}) we get, using (\ref{relentlow}) 
\begin{eqnarray}
q_{nL(l+m),\varepsilon } &\leq &\sum_{_{\substack{ -\frac{1}{nL}\log \lambda
\geq \underset{r\in A_{l,m,\varepsilon }^{c}}{\min }s_{(l+m)}(\widehat{\Psi }%
_{(r,l,m)})  \\ +\underset{r\in A_{l,m,\varepsilon }^{c}}{\min }s_{(l+m)}(%
\widehat{\Psi }_{(r,l,m)},\widehat{\Phi }_{(l,m)})-\varepsilon }}}\text{spec}%
_{\lambda }((\widehat{\Phi }_{(l,m)}^{(l+m)})^{\otimes nL})  \label{xxx} \\
&\leq &\sum_{_{\substack{ -\frac{1}{nL}\log \lambda \geq \underset{r\in
A_{l,m,\varepsilon }^{c}}{\min }s_{(l+m)}(\widehat{\Psi }_{(r,l,m)})  \\ %
+(l+m)s(\Psi ,\Phi )-(l+m+1)\varepsilon }}}\text{spec}_{\lambda }((\widehat{%
\Phi }_{(l,m)}^{(l+m)})^{\otimes nL}).  \notag
\end{eqnarray}%
Next, observe that by the subadditivity of the entropy we have for each $r$%
\begin{eqnarray*}
s_{(l+m)}(\widehat{\Psi }_{(r,l,m)}) &=&\lim_{k\rightarrow \infty }\frac{1}{k%
}S(\widehat{\Psi }_{(r,l,m)}^{(k(l+m))})=\lim_{k\rightarrow \infty }\frac{1}{%
k}S(\Psi _{(r,l+m)}\upharpoonright _{(\mathbf{I}^{\otimes l}\otimes \mathcal{%
A}^{\mathbb{\otimes }m})^{\mathbb{\otimes }k}}) \\
&\geq &\lim_{k\rightarrow \infty }\frac{1}{k}(S(\Psi
_{(r,l+m)}^{(k(l+m))})-kl\text{Tr}\mathbf{1}_{\mathcal{A}}) \\
&=&s_{(l+m)}(\Psi _{(r,l+m)})-l\text{Tr}\mathbf{1}_{\mathcal{A}}=(l+m)s(\Psi
)-l\text{Tr}\mathbf{1}_{\mathcal{A}},
\end{eqnarray*}%
and hence for sufficiently large $m$ we may continue the chain of
inequalities (\ref{xxx}): 
\begin{equation*}
q_{nL(l+m),\varepsilon }\leq \sum_{-\log \lambda \geq nL(l+m)(s(\Psi
)+s(\Psi ,\Phi )-3\varepsilon )}\text{spec}_{\lambda }((\widehat{\Phi }%
_{(l,m)}^{(l+m)})^{\otimes nL}).
\end{equation*}

From this we derive the following upper bound, being valid for $m$ large
enough (using the Araki-Lieb inequality in the fifth line) 
\begin{eqnarray*}
&&\widehat{\Phi }_{(l,m)}^{nL(l+m)}(q_{nL(l+m),\varepsilon }) \\
&\leq &e^{-nL(l+m)(s(\Psi )+s(\Psi ,\Phi )-3\varepsilon )}\text{Tr}%
q_{nL(l+m),\varepsilon } \\
&\leq &e^{-nL(l+m)(s(\Psi )+s(\Psi ,\Phi )-3\varepsilon )}\sum_{r\in
A_{l,m,\varepsilon }^{c}}\text{Tr}p_{r,nL,\varepsilon } \\
&\leq &e^{-nL(l+m)(s(\Psi )+s(\Psi ,\Phi )-3\varepsilon )}\sum_{r\in
A_{l,m,\varepsilon }^{c}}e^{nL(s_{(l+m)}(\widehat{\Psi }_{(r,l,m)})+%
\varepsilon )} \\
&\leq &e^{-nL(l+m)(s(\Psi )+s(\Psi ,\Phi )-3\varepsilon
)}(l+m)e^{nL(l+m)(s(\Psi )+2\varepsilon )} \\
&\leq &e^{-nL(l+m)(s(\Psi ,\Phi )-6\varepsilon )}.
\end{eqnarray*}

From $\ast $-mixing we get now the desired separation order 
\begin{eqnarray*}
\Phi ^{(nL(l+m))}(q_{nL(l+m),\varepsilon }) &\leq &e^{-nL(l+m)(s(\Psi ,\Phi
)-6\varepsilon )}\alpha ^{-nL} \\
&=&e^{-nL(l+m)(s(\Psi ,\Phi )+\frac{\log \alpha }{l+m}-6\varepsilon )} \\
&\leq &e^{-nL(l+m)(s(\Psi ,\Phi )-7\varepsilon )}
\end{eqnarray*}
(for $m$ large enough).

On the other hand, $\Psi $-typicality is guaranteed by 
\begin{eqnarray*}
\Psi ^{(nL(l+m))}(q_{nL(l+m),\varepsilon }) &=&\Psi ^{(nL(l+m))}\left(
\bigvee_{r\in A_{l,m,\varepsilon }^{c}}p_{r,nL,\varepsilon }\right) \\
&=&\frac{1}{l+m}\sum_{r^{\prime }=0}^{l+m-1}\Psi _{(r^{\prime
},l+m)}^{(nL(l+m))}\left( \bigvee_{r\in A_{l,m,\varepsilon
}^{c}}p_{r,nL,\varepsilon }\right) \\
&\geq &\frac{1}{l+m}\sum_{r^{\prime }\in A_{l,m,\varepsilon }^{c}}\Psi
_{(r^{\prime },l+m)}^{(nL(l+m))}\left( p_{r^{\prime },nL,\varepsilon }\right)
\\
&\geq &\frac{1}{l+m}\sum_{r^{\prime }\in A_{l,m,\varepsilon
}^{c}}(1-\varepsilon ),
\end{eqnarray*}
the last inequality being valid for large $n$. We may continue 
\begin{eqnarray*}
\Psi ^{(nL(l+m))}(q_{nL(l+m),\varepsilon }) &\geq &(1-\varepsilon )-\frac{1}{%
l+m}\sum_{r^{\prime }\in A_{l,m,\varepsilon }}(1-\varepsilon ) \\
&\geq &(1-\varepsilon )-\frac{\#A_{l,m,\varepsilon }}{l+m}\geq 1-2\varepsilon
\end{eqnarray*}
for $m$ large enough (by (\ref{denslem})).

Now (in the usual way) we may interpolate the $q_{nL(l+m),\varepsilon }$ in
order to define the projections $q_{n,\varepsilon }$ also for $n\in \mathbb{N%
}$ which are not multiples of $L(l+m)$. We derived the existence of a
sequence of projections being asymptotically $\varepsilon $-typical for $%
\Psi $ and fulfilling 
\begin{equation*}
\Phi ^{(n)}(q_{n,\varepsilon })\leq e^{-n(s(\Psi ,\Phi )-\varepsilon )}
\end{equation*}%
for large $n$. This proves that, for any $\alpha \in (0,1)$ the separation
exponent fulfils $\overline{\beta }_{\alpha }(\Psi ,\Phi )\leq -s(\Psi ,\Phi
)$ for finite $s(\Psi ,\Phi )$ .

5. Now assume $s(\Psi ,\Phi )=+\infty $. Observe that in that case the
estimates in (\ref{lowest1}) and hence (\ref{denslem}) are not valid. But (%
\ref{denslem}) becomes true if we replace the definition of $A_{l,m,\eta }$
most appropriately by 
\begin{equation*}
A_{l,m,\eta }:=\{r:0\leq r<l+m,\frac{1}{l+m}S(\widehat{\Psi }%
_{(r,l,m)}^{(l+m)},\widehat\Phi^{(m+l)}_{(m,l)} )<\eta ^{-1}\}.
\end{equation*}

In fact, choose $M$ large enough to ensure $S(\Psi ^{(M)},\Phi ^{(M)})>\eta
^{-1}M$ (we include the case $S(\Psi ^{(M)},\Phi ^{(M)})=+\infty $). Now we
have the ergodic decomposition 
\begin{equation}
\Psi =\frac{1}{M}\sum_{r=0}^{M-1}\Psi _{(r,M)}=\frac{1}{M}%
\sum_{r=0}^{M-1}\Psi _{(0,M)}\circ \tau ^{-r}  \label{ergdec}
\end{equation}
due to \cite{q-sm}. The states $\Psi _{(r,M)}=\Psi _{(0,M)}\circ \tau ^{-r}$
are $\tau ^{M}$-ergodic. In view of the (joint) convexity of the relative
entropy we conclude that at least one of the $r$ fulfils $S(\Psi
_{_{(r,M)}}^{(M)},\Phi ^{(M)})>\eta ^{-1}M$. We may assume without any loss
of generality that this is true for $r=0$, i.e. $S(\Psi
_{_{(0,M)}}^{(M)},\Phi ^{(M)})>\eta ^{-1}M$. The $\tau ^{M}$-ergodic state $%
\Psi _{(0,M)}$ again has an ergodic decomposition with respect to $\tau
^{2M} $%
\begin{equation*}
\Psi _{(0,M)}=\frac{1}{2}(\Psi ^{\prime }+\Psi ^{\prime }\circ \tau ^{-M})
\end{equation*}
and, applying once again the convexity argument we find that we may assume $%
S(\Psi ^{\prime (M)},\Phi ^{(M)})>\eta ^{-1}M$. $\Psi ^{\prime }$ is $\tau
^{2M}$-ergodic, and we obtain from (\ref{ergdec}) an ergodic decomposition
of $\Psi $ into $\tau ^{2M}$-ergodic states 
\begin{equation*}
\frac{1}{2M}\sum_{r=0}^{2M-1}\Psi ^{\prime }\circ \tau ^{-r},
\end{equation*}
hence we may assume without loss of generality that $\Psi _{(0,2M)}=\Psi
^{\prime }$. So we have 
\begin{equation*}
S(\Psi _{(0,2M)}^{(M)},\Phi ^{(M)})>\eta ^{-1}M
\end{equation*}
for $M$ large enough. This yields 
\begin{equation*}
S((\Psi _{(0,2M)}\circ \tau ^{-r})^{(\{r,r+1,...,r+M-1\})},\Phi ^{(M)})>\eta
^{-1}M\text{ \ for each }r
\end{equation*}
in view of the definition of $\tau $, i.e. (using the stationarity of $\Phi $%
) 
\begin{equation*}
S((\Psi _{(r,2M)})^{(\{r,r+1,...,r+M-1\})},\Phi
^{^{(\{r,r+1,...,r+M-1\})}})>\eta ^{-1}M.
\end{equation*}
In view of the monotonicity of the relative entropy we get now for $r\geq l$%
\begin{equation*}
S((\widehat{\Psi }_{(r,l,2M-l)})^{(2M)},\Phi \upharpoonright _{\mathbf{I}%
^{\otimes l}\otimes \mathcal{A}^{\mathbb{\otimes }m}})>\eta ^{-1}M.
\end{equation*}
So again (\ref{denslem}) is fulfilled for $M$ sufficiently large. We
conclude that asymptotically for the overwhelming part of the $r$ in $%
\{0,1,...,l+m-1\}$ the expression $\frac{1}{l+m}s_{(l+m)}(\widehat{\Psi }%
_{(r,l,m)},\widehat{\Phi }_{(l,m)})$ is arbitrary large (i.e. $>\frac{1}{2}%
\eta ^{-1}$ or even infinite) for large $m$. Now we may proceed essentially
as in 4., employing the results of \ \cite{igrai}. We find projections $%
p_{r,nL,\eta }$, separately for each $r$ in $A_{l,m,\eta }^{c}$, which
distinguish between $\widehat{\Psi }_{(r,l,m)}$ and $\widehat{\Phi }_{(l,m)}$
exponentially well at a rate at least$\frac{1}{3}\eta ^{-1}$, and we may
join these projections to find $\Psi $-typical projections $q_{nL(l+m),\eta }
$. This is possible due to the properties a) and b) above, where now b) is
modified to

b') The projections $p_{r,nL,\eta }$ are constructed as certain
sub-projections of the projection 
\begin{equation}
p^{(nL)}:=\sum_{\log \lambda \leq -nL\frac{1}{3}\eta ^{-1}}\text{spec}%
_{\lambda }((\widehat{\Phi }_{(l,m)}^{(l+m)})^{\otimes nL}),  \label{mod1}
\end{equation}

But we have to take into account that \cite{igrai} only treats the case of
finite relative entropy rate; hence those $r$, for which $S(\widehat{\Psi }%
_{(r,l,m)}^{(l+m)},\widehat{\Phi }_{(l,m)}^{(l+m)})=+\infty $, are still not
covered. For those $r$, simply choose 
\begin{equation*}
p_{r,nL,\eta }\equiv \text{spec}_{0}((\widehat{\Phi }_{(l,m)}^{(l+m)})^{%
\otimes nL}).
\end{equation*}
Obviously, this projection fulfils $(\widehat{\Phi }_{(l,m)}^{(l+m)})^{%
\otimes nL}($spec$_{0}((\widehat{\Phi }_{(l,m)}^{(l+m)})^{\otimes nL}))=0$,
and it is a sub-projection of $p^{(nL)}$. We still have to show that spec$%
_{0}((\widehat{\Phi }_{(l,m)}^{(l+m)})^{\otimes nL})$ is asymptotically
typical for each $\widehat{\Psi }_{(r,l,m)}$ with $S(\widehat{\Psi }%
_{(r,l,m)}^{(l+m)}, \widehat\Phi^{(m+l)}_{(m,l)})=+\infty $. In fact,
represent the density operator of $\Phi ^{(m)}$as 
\begin{equation*}
D_{\Phi ^{(m)}}=\sum_{j=1}^{K}\lambda _{j}w_{j}
\end{equation*}
where the $w_{j}$ are mutually orthogonal minimal projectors in $\mathcal{A}%
^{\mathbb{\otimes }m}$ fulfilling $\sum_{j}w_{j}=\mathbf{1}_{\mathcal{A}^{%
\mathbb{\otimes }m}}$ (and the $\lambda _{j}$ are the eigen-values of $%
D_{\Phi ^{(m)}}$ including $0$). Let $v_{j}:=\mathbf{1}_{\mathcal{A}^{%
\mathbb{\otimes }l}}\otimes w_{j}$. Observe that due to $S(\widehat{\Psi }%
_{(r,l,m)}^{(l+m)},\widehat{\Phi }_{(l,m)}^{(l+m)})=+\infty $ there is at
least one $\overline{j}$ with $\lambda _{\overline{j}}=0$ but $\widehat{\Psi 
}_{(r,l,m)}^{((l+m))}\left( v_{\overline{j}}\right) >0$. Now we have 
\begin{eqnarray}  \label{clcl}
\widehat{\Psi }_{(r,l,m)}^{(n(l+m))}(\text{spec}_{0}((\widehat{\Phi }%
_{(l,m)}^{(l+m)})^{\otimes n})) =\sum_{(j_{1},...,j_{n})\in N_{n}}\widehat{%
\Psi }_{(r,l,m)}^{(n(l+m))}\left( \bigotimes_{k=1}^{n}v_{j_{k}}\right) ,
\end{eqnarray}
where $N_{n}:=\{(j_{1},...,j_{n}):\prod_{k=1}^{n}\lambda
_{j_{k}}=0\}=((N_{1}^{c})^{n})^{c}$. Denote by $\mathcal{B}$ the abelian
sub-algebra of $\mathbf{I}^{\mathbb{\otimes }l}\otimes \mathcal{A}^{\mathbb{%
\otimes }m}$ generated by the set $\{v_{j}\}$. Then the quasi-local algebra $%
\mathcal{B}^{\mathbb{Z}}$ is an abelian sub-algebra of $(\mathbf{I}^{\mathbb{%
\otimes }l}\otimes \mathcal{A}^{\mathbb{\otimes }m})^{\mathbb{Z}}$ and the
restriction $P$ of $\widehat{\Psi }_{(r,l,m)}$ to this sub-algebra is a
classical ergodic process with $K$ symbols (Gelfand isomorphism and Riesz
representation theorem). This process fulfils $P^{(1)}(\{\overline{j}\})>0$.
We may continue the left-hand side in (\ref{clcl}) as follows 
\begin{eqnarray*}
\widehat{\Psi }_{(r,l,m)}^{(n(l+m))}(\text{spec}_{0}((\widehat{\Phi }
_{(l,m)}^{(l+m)})^{\otimes n})) &=&P^{(n)}\left( N_{n}\right) \\
&=&1-P^{(n)}((N_{1}^{c})^{n}) \\
&\geq &1-P^{(n)}((\{\overline{j}\}^{c})^{n}).
\end{eqnarray*}
Now $P^{(n)}((\{\overline{j}\}^{c})^{n})$ is the probability of all $n$%
-sequences of symbols where the symbol $\overline{j}$ does not appear at
all. This tends to zero, since by the individual ergodic theorem the a.s.
asymptotic frequency of the symbol $\overline{j}$ is $P^{(1)}(\{\overline{j}%
\})=\widehat{\Psi }_{(r,l,m)}^{((l+m))}\left( v_{\overline{j}}\right) >0$ by
assumption. Hence the conclusions of part 4. are valid in the case of
infinite relative entropy, too.
\end{proof}

\section{The stationary case\label{chstat}}

So far we formulated the Theorems \ref{stein-sanov} and \ref{classsan} for
sets of ergodic states $\Psi $ resp. processes $P$ to be optimally separated
from a reference state or process. These results can be easily extended to
the general stationary situation.

Any stationary state $\Psi \in \mathcal{S}_{\text{stat}}(\mathcal{A}^{%
\mathbb{Z}})$ can be represented as a mixture (ergodic decomposition) 
\begin{equation*}
\Psi =\int_{\mathcal{S}_{\text{erg}}(\mathcal{A}^{\mathbb{Z}})}\Xi \gamma
_{\Psi }(d\Xi )
\end{equation*}
of ergodic states ($\mathcal{S}_{\text{stat}}(\mathcal{A}^{\mathbb{Z}})$ is
a Choquet simplex, $\mathcal{S}_{\text{erg}}(\mathcal{A}^{\mathbb{Z}})$ is
the corresponding set of extremal points, $\gamma _{\Psi }$ is a probability
measure on the measurable space $[\mathcal{S}_{\text{erg}}(\mathcal{A}^{%
\mathbb{Z}}),\mathfrak{B}(\Upsilon _{\mathcal{A}^{\mathbb{Z}}})]$, with $%
\Upsilon _{\mathcal{A}^{\mathbb{Z}}}$ denoting the weak-$\ast $-topology and 
$\mathfrak{B}(\Upsilon _{\mathcal{A}^{\mathbb{Z}}})$ the corresponding Borel 
$\sigma $-field, cf. \cite{ruelle}). The measure $\gamma _{\Psi }$ is unique.

Now let $\Phi \in \mathcal{S}(\mathcal{A}^{\mathbb{Z}})$ be a state and $%
\Theta \subseteq \mathcal{S}_{\text{stat}}(\mathcal{A}^{\mathbb{Z}})$ with
the property that for any $\Psi \in \Theta $ the relative entropy rate $%
s(\Xi ,\Phi )$ exists for $\gamma _{\Psi }$-almost all $\Xi $. We define the
quantity 
\begin{eqnarray*}
\underline{s}(\Psi ,\Phi ):=\text{essinf}_{\gamma _{\Psi }(d\Xi )}s(\Xi
,\Phi ),
\end{eqnarray*}
and for $\Omega \subseteq \Theta $ the quantity 
\begin{eqnarray*}
\underline{s}(\Omega ,\Phi ):=\inf_{\Psi \in \Omega }\underline{s}(\Psi
,\Phi ).
\end{eqnarray*}

\begin{theorem}
\label{stat-sanov} Let $\Phi $ be a state on $\mathcal{A}^{\mathbb{Z}}$ and $%
\Theta \subseteq \mathcal{S}_{\text{stat}}(\mathcal{A}^{\mathbb{Z}})$ such
that for each $\Psi \in \Theta $ and $\gamma _{\Psi }$-almost all $\Xi $ the
pair $(\Xi ,\Phi )$ satisfies the HP-condition.

Then the quantity $\underline{s}(\Psi ,\Phi )\leq +\infty $ exists for each $%
\Psi \in \Theta $, and to each subset $\Omega \subseteq \Theta $ and any $%
\eta >0$ there exists a sequence $\{p_{n}\}_{n\in \mathbb{N}}$ of
projections $p_{n}\in \mathcal{A}^{(n)}$ with 
\begin{equation}
\lim_{n\rightarrow \infty }\Psi ^{(n)}(p_{n})=1,\qquad \text{for all }\Psi
\in \Omega  \label{stat1}
\end{equation}
and 
\begin{equation}
\underset{n\rightarrow \infty }{\text{\emph{limsup}}}\frac{1}{n}\log \Phi
^{(n)}(p_{n}) \leq -\underline{s}(\Omega ,\Phi )+\eta .  \label{stat2}
\end{equation}
if $\underline{s}(\Omega, \Phi) < \infty$, otherwise if $\underline{s}%
(\Omega, \Phi)= \infty$ 
\begin{equation}
\underset{n\rightarrow \infty }{\text{\emph{limsup}}}\frac{1}{n}\log \Phi
^{(n)}(p_{n}) \leq -\frac{1}{\eta}.  \label{stat2a}
\end{equation}
Moreover, for each sequence of projections $\{\widetilde{p}_{n}\}$
fulfilling (\ref{stat1}) we have 
\begin{equation}
\underset{n\rightarrow \infty }{\text{\emph{liminf}}}\frac{1}{n}\log \Phi
^{(n)}(\widetilde{p}_{n})\geq -\underline{s}(\Omega ,\Phi )\text{.}
\label{stat3}
\end{equation}
Hence $-\underline{s}(\Omega ,\Phi )$ is the lower limit of all achievable
separation exponents.
\end{theorem}

\begin{remark}
If $\Phi $ is stationary and, moreover, $\ast $-mixing, the assumption of
the Theorem \ref{stat-sanov} is fulfilled with $\Theta =\mathcal{S}_{\text{%
stat}}(\mathcal{A}^{\mathbb{Z}})$, according to section \ref{starmix}.
\end{remark}

\begin{proof}
Let $\tilde{\Omega }:=\{\Xi \in \mathcal{S}_{\text{erg}}(\mathcal{A}^{%
\mathbb{Z}}):(\Xi ,\Phi )$ satisfies the HP-condition and $s(\Xi ,\Phi )\geq 
\underline{s}(\Omega ,\Phi )\}$. The set $\tilde{\Omega }$ is weak-$\ast $%
-measurable since it can be represented by a countable application of unions
and intersections to local sets, defined via the measurable functions $%
S(\cdot ,\Phi ^{(n)})$ and $\beta _{\varepsilon ,n}(\cdot ,\Phi )$.

Let $p_{n}$ be chosen as in Theorem \ref{stein-sanov}, with $\Omega $ there
specified as $\widetilde{\Omega }$. Then (\ref{stat2}) or (\ref{stat2a}) are
trivially fulfilled. For any $\Psi \in \Theta $ we obtain by assumption 
\begin{eqnarray}
\Psi ^{(n)}(p_{n}) &=&\int_{\mathcal{S}_{\text{erg}}(\mathcal{A}^{\mathbb{Z}%
})}\Xi ^{(n)}(p_{n})\gamma _{\Psi }(d\Xi )  \label{ergrep1} \\
&=&\int_{\overline{\Omega }}\Xi ^{(n)}(p_{n})\gamma _{\Psi }(d\Xi ).  \notag
\end{eqnarray}%
Now for each $\Xi \in \overline{\Omega }$ the expression $\Xi
^{(n)}(p_{n})\in \lbrack 0,1]$ tends to $1$ by the choice of the projections 
$p_{n}$. Hence Lebesgue's theorem on dominated convergence guarantees (\ref%
{stat1}).

On the other hand, for each sequence of projections $\{\widetilde{p}_{n}\}$
fulfilling (\ref{stat1}) the identity (\ref{ergrep1}) (with $\widetilde{p}%
_{n}$ instead of $p_{n}$) proves that, for each $\Psi \in \Omega $, $\Xi
^{(n)}(\widetilde{p}_{n})$ tends to $1$ in $\gamma _{\Psi }$-probability as $%
n\rightarrow \infty $.\ By the definition of $\underline{s}(\Omega ,\Phi )$
to any $\eta >0$ we may choose $\Psi $ in such a way that $\underline{s}%
(\Psi ,\Phi )\leq \underline{s}(\Omega ,\Phi )+\eta $. We show that 
\begin{equation*}
\underset{n\rightarrow \infty }{\text{liminf}}\frac{1}{n}\log \Phi ^{(n)}(%
\widetilde{p}_{n})\geq -\underline{s}(\Psi ,\Phi ),
\end{equation*}
which implies (\ref{stat3}) since $\eta $ can be chosen arbitrarily small.
In fact, assume the existence of a sub-sequence $n^{\prime }$ such that 
\begin{equation}
\lim_{n^{\prime }}\frac{1}{n^{\prime }}\log \Phi ^{(n^{\prime })}(\widetilde{%
p}_{n^{\prime }})\leq -\underline{s}(\Psi ,\Phi )-\delta ,\delta >0.
\label{contra}
\end{equation}
Along that sub-sequence there is still convergence in $\gamma _{\Psi }$%
-probability of $\Xi ^{(n^{\prime })}(\widetilde{p}_{n^{\prime }})$ to $1$.
Since convergence in probability implies almost sure convergence of some
sub-sequence, we find another sub-sequence $n^{\prime \prime }$ of $%
n^{\prime }$ with $\lim_{n^{\prime \prime }}\Xi ^{(n^{\prime \prime })}(%
\widetilde{p}_{n^{\prime \prime }})=1$ holding $\gamma _{\Psi }$-almost
surely. Hence, in view of the definition of $\underline{s}(\Psi ,\Phi )$
there is some $\Xi _{0}\in \mathcal{S}_{\text{erg}}(\mathcal{A}^{\mathbb{Z}%
}) $ such that $(\Xi _{0},\Phi )$ fulfils the HP-condition, $s(\Xi _{0},\Phi
)<\underline{s}(\Psi ,\Phi )+\delta $, but

$\lim_{n^{\prime \prime }}\Xi _{0}^{(n^{\prime \prime })}(\widetilde{p}%
_{n^{\prime \prime }})=1$. Now Theorem \ref{stein-sanov}, applied to the
case $\Omega =\{\Xi _{0}\}$ implies 
\begin{equation*}
\underset{n^{\prime \prime }}{\text{liminf}}\frac{1}{n^{\prime \prime }}\log
\Phi ^{(n^{\prime \prime })}(\widetilde{p}_{n^{\prime \prime }})>-\underline{%
s}(\Psi ,\Phi )-\delta ,
\end{equation*}
which contradicts (\ref{contra}).
\end{proof}

The classical case immediately follows (with $\gamma _{P}$ denoting the
probability measure occurring in the ergodic decomposition of a stationary
process $P$ and $\underline{h}(P,Q):=$ess$\inf_{\gamma _{P}(dW)}h(W,Q)$,
supposed that $h(W,Q)$ exists $\gamma _{P}$-almost surely):

\begin{theorem}
\label{stat-sanov-class} Let $Q\in \mathcal{P}(A^{\mathbb{Z}})$ be a process
and $\Theta \subseteq \mathcal{P}_{\text{stat}}(A^{\mathbb{Z}})$. Assume
that for each $P\in \Theta $ and $\gamma _{P}$-almost all $W\in \mathcal{P}_{%
\text{stat}}(A^{\mathbb{Z}})$ the relative entropy rate $h(W,Q)$ exists and $%
\bar{\beta}_{\varepsilon }(W,Q)\leq -h(W,Q)$ for all $\varepsilon \in (0,1)$.

Then the quantity $\underline{h}(P,Q)\leq +\infty $ exists for all $P\in
\Theta $, and to each subset $\Omega \subseteq \Theta $ and any $\eta >0$
there exists a sequence $\{M_{n}\}_{n\in \mathbb{N}}$ of subsets $%
M_{n}\subseteq A^{n}$ with 
\begin{equation}
\lim_{n\rightarrow \infty }P^{(n)}(M_{n})=1,\qquad \text{for all }P\in \Omega
\label{cccl1}
\end{equation}
and 
\begin{equation*}
\underset{n\rightarrow \infty }{\text{\emph{limsup}}}\frac{1}{n}\log
Q^{(n)}(M_{n}) \leq -\underline{h}(\Omega ,Q)+\eta
\end{equation*}
if $\underline{h}(\Omega, Q)< \infty$, otherwise if $\underline{h}(\Omega,
Q)=\infty$ 
\begin{equation*}
\underset{n\rightarrow \infty }{\text{\emph{limsup}}}\frac{1}{n}\log
Q^{(n)}(M_{n}) \leq -\frac{1}{\eta}.
\end{equation*}
Moreover, for each sequence of subsets $\{\widetilde{M}_{n}\}$ fulfilling (%
\ref{cccl1}) we have 
\begin{equation*}
\underset{n\rightarrow \infty }{\text{\emph{liminf}}}\frac{1}{n}\log Q^{(n)}(%
\widetilde{M}_{n}) \geq -\underline{h}(\Omega ,Q)\text{.}
\end{equation*}
Hence $-\underline{h}(\Omega ,Q)$ is the lower limit of all achievable
separation exponents.
\end{theorem}

\begin{remark}
If $Q$ is stationary and, moreover, $\ast $-mixing, the assumption of the
Theorem is fulfilled with $\Theta =\mathcal{P}_{\text{stat}}(A^{\mathbb{Z}})$%
, according to section \ref{starmix}.
\end{remark}

\section{The quantum Shannon-McMillan theorem for stationary states and
other corollaries\label{chcoro}}

As announced in the introduction, several earlier results on typical
subspaces resp. subsets are contained in Theorem \ref{stat-sanov} in a
version extended to the stationary case. We emphasize that the initial
versions of quantum Shannon-McMillan theorem, Kaltchenko-Yang universality
and the quantum Stein's lemma were important ingredients in our proof. Also,
it should be mentioned that it is not difficult to prove the stationary case
of the quantum Shannon-McMillan theorem directly from the Kaltchenko-Yang
result, without using quantum Sanov's theorem.

\begin{corollary}
\label{statqsm}(Quantum Shannon-McMillan theorem for stationary states)

Let $\Psi \in \mathcal{S}_{\text{stat}}(\mathcal{A}^{\mathbb{Z}})$ be a
stationary state and $\Psi =\int_{\mathcal{S}_{\text{erg}}(\mathcal{A}^{%
\mathbb{Z}})}\Xi \gamma _{\Psi }(d\Xi )$ be its ergodic decomposition. Then
there exists a sequence $\{p_{n}\}$ of projections in $\mathcal{A}^{(n)}$,
respectively such that

\begin{itemize}
\item $\lim_{n\rightarrow \infty }\Psi ^{(n)}(p_{n})=1$ (typicality)

\item $\lim_{n\rightarrow \infty }\frac{1}{n}$\emph{Tr}$p_{n}=$\emph{esssup}$%
_{\gamma _{\Psi }(d\Xi )}s(\Xi ):=\overline{s}(\Psi )$ (maximal ergodic
entropy rate).
\end{itemize}

For any sequence $\widetilde{p}_{n}$ with $\lim_{n\rightarrow \infty }\Psi
^{(n)}(\widetilde{p}_{n})=1$ we have 
\begin{eqnarray*}
\underset{n\rightarrow \infty }{\liminf}\frac{1}{n} \text{Tr} \widetilde{p}%
_{n}\geq \overline{s}(\Psi ) \text{ (optimality).}
\end{eqnarray*}
\end{corollary}

\begin{remark}
We emphasize that the AEP does not hold in the stationary case.

Also, observe that the relevant notion in the stationary case is \emph{not}
the von Neumann entropy rate $s(\Psi )$ of the state $\Psi $ being the
average of the entropy rates of the ergodic components of $\Psi $, but their
essential supremum $\overline{s}(\Psi )$.
\end{remark}

\begin{proof}
Let $\Phi $ be the tracial state in $\mathcal{S}(\mathcal{A}^{\mathbb{Z}})$.
It is $\ast $-mixing (even iid). Apply Theorem \ref{stat-sanov} with $\Omega
=\{\Psi \}$. This yields a sequence $\{p_{n}^{\eta }\}$ of $\Psi $-typical
projections with 
\begin{equation*}
\overline{s}(\Psi )\leq \underset{n\rightarrow \infty }{\text{liminf}}\frac{1%
}{n}\text{Tr}p_{n}^{\eta }\leq \underset{n\rightarrow \infty }{\text{limsup}}%
\frac{1}{n}\text{Tr}p_{n}^{\eta }\leq \overline{s}(\Psi )+\eta
\end{equation*}
for any $\eta >0$. Now the assertion of the Theorem easily follows, since $%
\Omega $ is a finite set.
\end{proof}

The next corollary extends the universality result of \cite{ky} to
stationary states:

\begin{corollary}
(Kaltchenko-Yang universality theorem for stationary states)

Let $\Omega _{s}:=\{\Psi \in \mathcal{S}_{\text{stat}}(\mathcal{A}^{\mathbb{Z%
}}):\overline{s}(\Psi )<s\}$. Then there exists a sequence $\{p_{n}\}$ of
projections in $\mathcal{A}^{(n)}$, respectively such that

\begin{itemize}
\item $\lim_{n\rightarrow \infty }\Psi ^{(n)}(p_{n})=1$ for each $\Psi \in $ 
$\Omega _{s}$ (typicality)

\item $\lim_{n\rightarrow \infty }\frac{1}{n}$\emph{Tr}$p_{n}=s$ (maximal
ergodic entropy rate).
\end{itemize}

For any sequence $\widetilde{p}_{n}$ with $\lim_{n\rightarrow \infty }\Psi
^{(n)}(\widetilde{p}_{n})=1,\Psi \in $ $\Omega _{s}$, we have 
\begin{equation*}
\underset{n\rightarrow \infty }{\text{\emph{liminf}}}\frac{1}{n}{\text{\emph{%
Tr}}}\widetilde{p}_{n}\geq s\text{ \ \ \ (optimality).}
\end{equation*}
\end{corollary}

\begin{proof}
Let $\Phi $ again be the tracial state in $\mathcal{S}(\mathcal{A}^{\mathbb{Z%
}})$. Apply Theorem \ref{stat-sanov} in a similar way as in the proof of
Corollary \ref{statqsm} to the sets $\Omega _{s-\eta },\eta >0$.
\end{proof}

\begin{remark}
Observe that the condition $\overline{s}(\Psi )<s$ defining $\Omega _{s}$
cannot be replaced by $\overline{s}(\Psi )\leq s$.
\end{remark}

Finally, quantum Stein's lemma \cite{igrai} is extended to the case where
the null hypothesis state $\Psi $ is only assumed stationary, the reference
state $\Phi $\ fulfills the HP-condition with respect to almost all ergodic
components of $\Psi $ (and the relative entropy rate $s(\Psi,\Phi)$ may be
infinite):

\begin{corollary}
(Stein's lemma for stationary states)

Let $\Phi \in \mathcal{S}(\mathcal{A}^{\mathbb{Z}})$ and $\Psi \in \mathcal{S%
}_{\text{stat}}(\mathcal{A}^{\mathbb{Z}})$ such that for $\gamma _{\Psi }$%
-almost all $\Xi $ the HP-condition \ is fulfilled for $(\Xi,\Phi)$. Then
there exists a sequence $\{p_{n}\}$ of projections with

\begin{itemize}
\item $\lim_{n\rightarrow \infty }\Psi ^{(n)}(p_{n})=1$ (typicality)

\item $\lim_{n\rightarrow \infty }\frac{1}{n}\log \Phi ^{(n)}(p_{n})= -%
\underline{s}(\Psi,\Phi)$ (achievability of the separation exponent $-%
\underline{s}(\Psi,\Phi)$).
\end{itemize}

For any sequence $\widetilde{p}_{n}$ with $\lim_{n\rightarrow \infty }\Psi
^{(n)}(\widetilde{p}_{n})=1$ we have 
\begin{eqnarray*}
\underset{n\rightarrow \infty }{\text{{\emph{liminf}}}}\frac{1}{n}\log \Phi
^{(n)}(\widetilde{p}_{n})\geq -\underline{s}(\Psi ,\Phi ) \text{ (optimality)%
}.
\end{eqnarray*}
\end{corollary}

\begin{remark}
Note that the relative AEP does not hold in Stein's lemma in the stationary
case.

Again, the relevant quantity in the stationary situation is not the average
relative entropy rate $s(\Psi ,\Phi )$, but the essential infimum \underline{%
$s$}$(\Psi ,\Phi )$.
\end{remark}

\begin{proof}
With $\Omega $ consisting of a single state $\Psi $ only, we may proceed in
the same way as in the proof of Corollary \ref{statqsm}.
\end{proof}

\section{Conclusions}

The paper is devoted to a generalization of Sanov's theorem from the iid
classical situation to the correlated case and, moreover, to the quantum
setting. In the present form, the main result comprises and extends several
earlier assertions including the (quantum and classical) Shannon-McMillan
theorem, Stein's Lemma (with relative AEP), Kaltchenko and Yang's
universality and, of course, a version of Sanov's theorem itself. It is a
continuation of \cite{sanov-i.i.d}, where the uncorrelated case is
considered. It has to be pointed out again (see \cite{sanov-i.i.d}), that
any attempt to formulate a quantisized version of Sanov's result has to face
the problem that the very notion of a trajectory and its empirical
distribution is problematic in quantum mechanics. Sanov's classical theorem
claims that for an iid process with marginal $Q$ the probability to produce
a trajectory with the empirical distribution belonging to some set $\Omega $
of probability measures is (in general) exponentially small. The
corresponding rate is specified as the minimal relative entropy between $Q$
and the distributions in $\Omega $. In the interesting case the measure $Q$
is of course \emph{not }an element of $\Omega $ or its topological closure.
So it is a large deviation result: the typical behaviour of $Q$-trajectories
is to have an empirical distribution close to $Q$. Whatever one tries to
adopt as a quantum substitute for the empirical distribution, the natural
choice in the case of a tensor product of vector states $v\otimes v\otimes
...\otimes v$ should be $v$ itself. This leads into the problem that for a
reference vector state $w^{\otimes n}$ the probability of measuring an
'empirical state' $v$ is at least $\text{Tr}P_{w^{\otimes n}}P_{v^{\otimes
n}}=|\langle w|v\rangle |^{2n}$, while the relative entropy of $v$ wrt $w$
is infinite, which would imply a super-exponential decay; for a more
detailed exposition see \cite{sanov-i.i.d}. In this situation it proves
useful to look at Sanov's theorem as an assertion about the likelihood of
observing the classical iid process given by $Q$ far from its original
support in the vicinity of the supports of \emph{other} iid processes. The
most natural choice for 'typical support' in the classical framework is the
set of trajectories with empirical distribution close to the given
probability distribution, since according to the individual ergodic theorem
the empirical distribution tends to $Q$ with probability one. So Sanov's
theorem in its original form says that the probability of observing the
trajectory in the typical support of other distributions, \emph{concretely
specified }by means of the corresponding empirical distributions, vanishes
at a rate given by the minimum relative entropy. It is of course completely
legitimate to insist on the point of view, that a quantum Sanov's theorem
should be about empirical distributions, too (see \cite{nagaokahayashi},
Remark 4, see also an attempt to formulate a quantum (iid) Sanov theorem
made in Segre's Ph.D. thesis \cite{segre} (2004), Conjecture 7.3.1.). But,
as explained, then one loses the relation to the established form of quantum
relative entropy (Umegaki's relative entropy). We chose to 'sacrifice'
empirical distributions in our approach but nonetheless calling it a version
of Sanov's theorem: in the classical case the relative entropy is not only
the rate of separation when empirical distributions as specifying typical
sets are considered. It has a clear operational meaning as the \emph{optimal
separation rate}, whatever one considers as typical support in the sense
that the probability goes to one. This perception of Sanov's theorem,
closely connected with the statistical hypothesis testing aspect, appears to
be natural. It allows useful generalizations to the correlated and quantum
cases.

\emph{Acknowledgements. }This work was supported by DFG grants "Entropy and
coding of large quantum systems" and by the Max-Planck Institute for
Mathematics in the Sciences, Leipzig. Tyll Kr\"{u}ger, Rainer
Siegmund-Schultze and Arleta Szko\l a are particularly grateful to Nihat Ay
for his constant encouragement during the preparation of the manuscript.

\end{document}